\documentclass[a4paper,10pt]{article}
\usepackage[utf8x]{inputenc}
\usepackage{tracefnt,amsmath,tabu,array}
\usepackage{amssymb,graphicx,setspace,amsfonts,amsbsy}
\usepackage{pifont,latexsym,ifthen,amsthm,rotating,calc,textcase,booktabs,cancel,slashed}

\newtheorem{theorem}{Theorem}[section]
\newtheorem{lemma}[theorem]{Lemma}
\newtheorem{corollary}[theorem]{Corollary}
\newtheorem{proposition}[theorem]{Proposition}
\newtheorem{definition}[theorem]{Definition}

\newtheorem{remark}[theorem]{Remark}
\newcommand{\filledbox}{\leavevmode
  \hbox to.77778em{%
  \hfil\vbox to.675em{\hrule width.6em height.6em}\hfil}}

\newcommand{\Rm}{{\mathbb R}}

\newcommand{\eps}{\varepsilon}

\begin{document}
\tabulinesep=1.0mm
\title{Inward/outward Energy Theory of Wave Equation in Higher Dimensions\footnote{MSC classes: 35L05, 35L71; This work is supported by National Natural Science Foundation of China Programs 11601374, 11771325}}

\author{Ruipeng Shen\\
Centre for Applied Mathematics\\
Tianjin University\\
Tianjin, China}

\maketitle

\begin{abstract}
  We consider the semi-linear, defocusing wave equation $\partial_t^2 u - \Delta u = -|u|^{p-1} u$ in $\Rm^d$ with $1+4/(d-1)\leq p < 1+4/(d-2)$. We generalize the inward/outward energy theory and weighted Morawetz estimates in 3D to higher dimensions. As an application we show the scattering of solutions if the energy of initial data decays at a certain rate as $|x| \rightarrow \infty$. 
\end{abstract}

\section{Introduction}

\subsection{Background}
In this work we consider defocusing wave equation in energy subcritical case $p < 1+4/(d-2)$
\[
 \left\{\begin{array}{ll} \partial_t^2 u - \Delta u = - |u|^{p-1}u, & (x,t) \in \Rm^d \times \Rm; \\
 u(\cdot, 0) = u_0; & \\
 u_t (\cdot,0) = u_1. & \end{array}\right.\quad (CP1)
\]
\paragraph{Critical Sobolev space} This wave equation is invariant under the natural rescaling transformation, i.e. if $u$ is a solution to (CP1), then the function
\[
 u_\lambda (x,t) = \lambda^{-2/(p-1)}u(x/\lambda,t/\lambda) 
\]
is another solution to (CP1). In addition we have the identity
\begin{align*}
\left\|(u_\lambda(\cdot,\lambda t'), \partial_t u_\lambda (\cdot,\lambda t'))\right\|_{\dot{H}^{s_p} \times \dot{H}^{s_p-1}} = \left\|(u(\cdot,t'), \partial_t u(\cdot,t'))\right\|_{\dot{H}^{s_p} \times \dot{H}^{s_p-1}},
\end{align*}
if we choose $s_p = d/2 - 2/(p-1)$. Thus the space $\dot{H}^{s_p}\times \dot{H}^{s_p-1}$ is called the critical Sobolev space of this equation. In particular, if $p=p_e(d) \doteq 1+\frac{4}{d-2}$, then $s_p=1$, the critical Sobolev space is exactly the energy space $\dot{H}^1 \times L^2$. We call this case the energy critical case. If $p= p_c(d) \doteq 1+\frac{4}{d-1}$, then $s_p = 1/2$, the critical Sobolev space is $\dot{H}^{1/2}\times \dot{H}^{-1/2}$. This case is usually called the conformal case. 

\paragraph{Local theory} The existence and uniqueness of solutions to (CP1) has been known for many years. Please see Kapitanski \cite{loc1} and Lindblad-Sogge \cite{ls}, for example. The idea is to combine suitable Strichartz estimates with a fixed-point argument. One may refer to Ginibre-Velo \cite{strichartz} for an almost complete version of Strichartz estimates and Keel-Tao \cite{endpointStrichartz} for a few endpoint cases. 

\paragraph{Scattering} The global behaviour of solutions has also been extensively studied in recent years. In general, we conjecture that all solutions to defocusing wave equation with initial data in the critical Sobolev spaces always scatter, i.e. a solution to defocusing nonlinear equation becomes more and more like a solution to the linear wave equation as $t$ tends to infinity. This was first proved in the energy critical case in 1990's. Please see Grillakis \cite{mg1} for dimension 3 and Grillakis \cite{mg2}, Shatah-Struwe \cite{ss1,ss2} for higher dimensions. If $p$ is energy supercritical ($p>p_e(d)$) or energy subcritical ($p<p_e(d)$), this conjecture becomes more difficult and stays to be an open problem, in spite of some progress, most in the 3-dimensional case.
\begin{itemize}
 \item There are many conditional results that prove the scattering if the critical Sobolev norm of solutions is uniformly bounded in the whole maximal lifespan, for different dimensions and ranges of $p$. For example, one may refer to Duyckaerts et al. \cite{dkm2}, Kenig-Merle \cite{km}, Killip-Visan \cite{kv2} (dimension 3), Killip-Visan \cite{kv3} (all dimensions) for energy supercritical case and Dodson-Lawrie \cite{cubic3dwave}, Shen \cite{shen2}, Dodson et al. \cite{nonradial3p5} (dimension 3) for energy subcritical case. All these works use the compactness-rigidity argument, which was first introduce by Keng-Merle \cite{kenig, kenig1} to study the global behaviour of solution to focusing, energy critical wave and Schr\"{o}dinger equations. This argument works in the focusing case as well.
 \item There are also many works proving the scattering of solutions under additional assumptions on the initial data. For example, if $d \geq 3$ and $p \in [p_c(d),p_e(d))$,  one may apply the conformal conservation law to prove the scattering of solutions if the initial data $(u_0,u_1)$ satisfy
 \[
  \int_{\Rm^d} \left[(|x|^2+1) (|\nabla u_0 (x)|^2 + |u_1(x)|^2) + |u_0(x)|^2 \right] dx < \infty,
 \]
 as shown in Ginibre-Velo \cite{conformal2} and Hidano \cite{conformal}. For another example, Yang \cite{yang1} proves the scattering under a weaker assumption on initial data
 \[
  \int_{\Rm^d} (1+|x|)^\gamma \left(\frac{1}{2}|\nabla u_0(x)|^2 + \frac{1}{2}|u_1(x)|^2 + \frac{1}{p+1}|u_0(x)|^{p+1}\right) < +\infty,
 \]
with $p$ and $\gamma$ satisfying 
\begin{align*}
 &\frac{1+\sqrt{d^2+4d-4}}{d-1} < p < p_e(d),& &\gamma> \max\left\{\frac{4}{p-1}-d+2,1\right\};&
\end{align*}
by considering the energy momentum tensor and its associated current. Recently Dodson \cite{claim1} proves the global existence and scattering of solutions in the conformal case of dimension 3 ($d=p=3$) under the assumption that the initial data $(u_0,u_1) \in \dot{H}^{1/2} \times \dot{H}^{-1/2}$ are radial.  
\end{itemize} 

\paragraph{Inward/outward energy theory} The author introduced an inward/outward energy theory on defocusing 3-dimensional wave equation in recent works \cite{shenenergy, sheninward, shen3dnonradial}. The exponent is assumed to be super-conformal but sub-critical, i.e. $p \in [3,5)$. The idea is to consider the inward/outward energies
\begin{align*}
 E_-(t) = \int_{\Rm^3} \left(\frac{1}{4}\left|\frac{x}{|x|}\cdot \nabla u + \frac{u}{|x|} + u_t \right|^2 + \frac{1}{4}|\slashed{\nabla} u|^2 + \frac{1}{2(p+1)}|u|^{p+1} \right) dx;\\
 E_+(t) = \int_{\Rm^3} \left(\frac{1}{4}\left|\frac{x}{|x|}\cdot \nabla u + \frac{u}{|x|} - u_t \right|^2 + \frac{1}{4}|\slashed{\nabla} u|^2 + \frac{1}{2(p+1)}|u|^{p+1} \right) dx.
\end{align*} and their fluxes through certain hyper-surfaces. Here the notation $\slashed{\nabla}$ represents the covariant derivative on the sphere centred at the origin. Thus $|\slashed{\nabla} u|^2 = |\nabla_x u|^2 - |\partial_r u|^2$. The consequences of this theory include
\begin{itemize} 
 \item The asymptotic behaviour of inward/outward energies
 \begin{align*}
  &\lim_{t \rightarrow +\infty} E_-(t) = 0;& &\lim_{t \rightarrow -\infty} E_+(t) = 0;& &\Rightarrow \lim_{t\rightarrow \pm \infty}\int_{\Rm^3} \left(|\slashed{\nabla} u|^2 + |u|^{p+1}\right) dx = 0.&
 \end{align*}
 \item Weighted Morawetz estimates. Given $0<\kappa<1$, we have
 \begin{align*}
  \int_{-\infty}^{\infty} \int_{\Rm^3} (|t|+|x|)^\kappa & \left(\frac{|\slashed{\nabla} u|^2}{|x|} + \frac{|u|^{p+1}}{|x|}\right) dx dt\\
  & \lesssim  E_\kappa(u_0,u_1) \doteq \int_{\Rm^3} |x|^\kappa \left(\frac{|\nabla u_0|^2}{2} + \frac{|u_1|^2}{2} + \frac{|u_0|^{p+1}}{p+1}\right) dx
 \end{align*}
 This immediately gives us the scattering of solutions if the initial data $(u_0,u_1)$ satisfy $E_\kappa(u_0,u_1)<+\infty$ for some $\kappa > \frac{5-p}{2}$. This assumption on the decay rate of initial data is weaker than all previously known results without radial assumptions. 
 \item If the initial data are radial, we may prove the scattering of solutions in the energy space $\dot{H}^1 \times L^2$ as long as $E_\kappa(u_0,u_1)<+\infty$ for some $\kappa \geq \frac{5-p}{p+1}$. Please note that this assumption on the decay rate is so weak that it can not guarantee $(u_0,u_1) \in \dot{H}^{s_p} \times \dot{H}^{s_p-1}$. Therefore we discover a scattering phenomenon that can not be covered by any previously known scattering theory. 
\end{itemize}

\paragraph{Main topic of this paper} In this work we generalize our inward/outward energy theory to higher dimensions. Our conclusion is similar to the 3-dimensional case. The general idea and major steps of argument remain the same as well. Nevertheless, we still introduce a few new ingredients in our argument. 
\begin{itemize}
 \item We give an almost complete version of Morawetz estimates. It is well known that if $d \geq 3$, the solution to (CP1) with a finite energy $E$ satisfies 
 \[
  \int_{-\infty}^\infty \int_{\Rm^d} \frac{|u(x,t)|^{p+1}}{|x|} dx dt \lesssim E.
 \]
Please see, for instance, Perthame-Vega \cite{benoit}. The integral of $|\slashed{\nabla} u(x,t)|^2/|x|$ in the whole space-time is also dominated by the energy $E$. This fact has been explained in \cite{shen3dnonradial} and independently proved in Yang \cite{yang1}. In subsection \ref{sec: Morawetz subsection} below we show that one more term, i.e. the integral of $|u(x,t)|^2/|x|^3$, can be inserted to the left hand of Morawetz inequality if $d \geq 4$. This plays a crucial role in the inward/outward energy theory. In fact, the value of the Morawetz integral 
\[
 \int_{t_1}^{t_2} \int_{\Rm^d} \left(\frac{1}{2} \cdot \frac{|\slashed{\nabla} u|^2}{|x|} + \frac{(d-1)(p-1)}{4(p+1)}\cdot \frac{|u|^{p+1}}{|x|} + \frac{(d-1)(d-3)}{8} \cdot \frac{|u|^2}{|x|^3}\right) dx dt 
\]
is exactly the amount of inward energy transformed to outward energy in the given period of time, as shown in Subsection \ref{statement eff}. The inward/outward energy theory also explains why our Morawetz estimates are ``almost complete''. Please see Proposition \ref{rediscover of Morawetz}. 
 \item We give a new proof of the energy flux formula. We use the operators $\mathbf{L}_{ij} = x_i \partial_{x_j} - x_j \partial_{x_i}$ rather than specific spherical coordinates, in order to simplify the calculation in higher dimensions. This also helps to give an energy flux formula of inward/outward energy for a non-radial region, as indicated in Proposition \ref{nonradial regions}.
 \item Second derivatives are involved in the calculation when we prove Morawetz estimates or energy flux formula of inward/outward energy. Thus we need to apply approximation techniques if the solution is not sufficiently smooth. This necessity becomes more obvious in higher dimensions. Because the regularity of the function $F(u) =-|u|^{p-1}u$ at $u=0$ depends on the exponent $p$. A higher dimension means smaller exponents $p$ and potentially less regular nonlinear term. In Remark \ref{smooth approximation Morawetz} and \ref{smooth approximation flux formula} below we give details of a useful smooth approximation technique, which works under a technical assumption
 \[
  F(u) = -|u|^{p-1} u \in L_{loc}^1 L^2 (\Rm \times \Rm^d) \quad \Leftrightarrow \quad u \in L_{loc}^p L^{2p} (\Rm \times \Rm^d).
 \]
 This is true for all finite-energy solutions as long as $(d,p)$ satisfies
 \[
  \left\{\begin{array}{ll} 1+\frac{4}{d-1}\leq p <  1+ \frac{4}{d-2}, & \hbox{if}\; 3\leq d\leq 6;\cr
  1+\frac{4}{d-1} \leq p \leq  1+\frac{3}{d-3}, & \hbox{if}\; 7\leq d\leq 9. 
  \end{array}\right. \quad (A1)
 \]
 For other pairs $(d,p)$, however, the same argument fails and we are not able to verify this technical assumption. Both these two situations are discussed in Remark \ref{local L1L2}. As a result we always assume that $(d,p)$ satisfies (A1) in this work. Nevertheless, the author believes that the inward/outward energy theory might probably be true as well for all $d \geq 3$ and $p \in [p_c(d),p_e(d))$, if we manage to develop a finer version of smooth approximation techniques. 
\end{itemize}

\subsection{Main Results}

We start by introducing the following notations. Here $r = |x|$ is the radius in $\Rm^d$.

\begin{definition}
 We define 
 \begin{align*}
  \mathbf{L} u & = r^{-\frac{d-1}{2}} \cdot \frac{\partial}{\partial_r}\left(r^{\frac{d-1}{2}} u\right) = \frac{x}{|x|}\cdot \nabla u + \frac{d-1}{2}\cdot \frac{u}{|x|};\\
  \mathbf{L}_+ u & = r^{-\frac{d-1}{2}} \cdot \frac{\partial}{\partial_r}\left(r^{\frac{d-1}{2}} u\right) + u_t = \frac{x}{|x|}\cdot \nabla u + \frac{d-1}{2}\cdot \frac{u}{|x|} + u_t;\\
  \mathbf{L}_- u & = r^{-\frac{d-1}{2}} \cdot \frac{\partial}{\partial_r}\left(r^{\frac{d-1}{2}} u\right) - u_t = \frac{x}{|x|}\cdot \nabla u + \frac{d-1}{2}\cdot \frac{u}{|x|} - u_t.
 \end{align*}
 When the initial data are involved we also use the notation 
 \begin{align*}
  \mathbf{L}_+ (u_0,u_1) & =  \frac{x}{|x|}\cdot \nabla u_0 + \frac{d-1}{2}\cdot \frac{u_0}{|x|} + u_1;\\
  \mathbf{L}_- (u_0,u_1) & =  \frac{x}{|x|}\cdot \nabla u_0 + \frac{d-1}{2}\cdot \frac{u_0}{|x|} - u_1;
 \end{align*}
\end{definition}

\begin{definition}
Let $u$ be a solution to (CP1) with a finite energy. Then we always have $(u(\cdot,t), u_t (\cdot,t)) \in (\dot{H}^1 \cap L^{p+1}) \times L^2$. We define inward/outward energy at time $t$ by
  \begin{align*}
  E_-(t) & = \int_{\Rm^d} \left(\frac{1}{4} |\mathbf{L}_+ u(x,t)|^2 + \frac{\lambda_d}{4}\cdot \frac{|u(x,t)|^2}{|x|^2} + \frac{1}{4}|\slashed{\nabla} u(x,t)|^2 + \frac{1}{2(p+1)}|u(x,t)|^{p+1} \right) dx;\\
  E_+(t) & = \int_{\Rm^d} \left(\frac{1}{4} |\mathbf{L}_- u(x,t)|^2 + \frac{\lambda_d}{4}\cdot \frac{|u(x,t)|^2}{|x|^2} + \frac{1}{4}|\slashed{\nabla} u(x,t)|^2 + \frac{1}{2(p+1)}|u(x,t)|^{p+1} \right) dx.
 \end{align*}
Here the constant $\lambda_d \doteq (d-1)(d-3)/4\geq 0$ depends on the dimension $d \geq 3$. We may also consider the inward/outward energy in a region $\Sigma \subset \Rm^d$ at time $t$ 
 \begin{align*}
  E_-(t;\Sigma) & = \int_{\Sigma} \left(\frac{1}{4} |\mathbf{L}_+ u(x,t)|^2 + \frac{\lambda_d}{4}\cdot \frac{|u(x,t)|^2}{|x|^2} + \frac{1}{4}|\slashed{\nabla} u(x,t)|^2 + \frac{1}{2(p+1)}|u(x,t)|^{p+1} \right) dx;\\
  E_+(t;\Sigma) & = \int_{\Sigma} \left(\frac{1}{4} |\mathbf{L}_- u(x,t)|^2 + \frac{\lambda_d}{4}\cdot \frac{|u(x,t)|^2}{|x|^2} + \frac{1}{4}|\slashed{\nabla} u(x,t)|^2 + \frac{1}{2(p+1)}|u(x,t)|^{p+1} \right) dx.
 \end{align*}
\end{definition}
\begin{remark}
 We have $E_-(t) +E_+(t) = E$ for all time $t$, according to Lemma \ref{identity of w u energy}. However, $if \Sigma \neq \Rm^d$, in general we may have
 \[
  E_-(\Sigma;t) + E_+(\Sigma;t) \neq \int_{\Sigma} \left(\frac{1}{2}|\nabla u(x,t)|^2 + \frac{1}{2}|u_t(x,t)|^2 + \frac{1}{p+1} |u(x,t)|^{p+1}\right) dx.
 \]
\end{remark}

\paragraph{Statements of main theorems} Now we give the main results of this paper.

\begin{theorem} \label{main 1}
Assume that $(d,p)$ satisfies (A1). Let $u$ be a solution to (CP1) with a finite energy $E$. Then we have the following asymptotic behaviour regarding the energy of $u$
\begin{align*}
& \lim_{t \rightarrow \pm \infty} E_{\mp} (t) = 0 \; \Rightarrow \; \lim_{t \rightarrow \pm \infty} \int_{\Rm^d} \left(\frac{\lambda_d |u(x,t)|^2}{4|x|^2} + \frac{|\slashed{\nabla} u(x,t)|^2}{4} + \frac{|u(x,t)|^{p+1}}{2(p+1)} \right) dx = 0.\\
& \lim_{t \rightarrow \pm \infty} \int_{|x|<c|t|} \left(\frac{1}{2}|\nabla u(x,t)|^2 + \frac{1}{2}|u_t(x,t)|^2 + \frac{1}{p+1} |u(x,t)|^{p+1}\right) dx = 0, \quad c\in (0,1).
\end{align*}
\end{theorem}

\begin{theorem}\label{main 2}
 Assume that $(d,p)$ satisfies assumption (A1) and $0<\kappa<1$. Let $(u_0,u_1)$ be initial data satisfying
 \[
  E_\kappa(u_0,u_1) \doteq \int_{\Rm^d} (1+|x|^{\kappa})\left[\frac{1}{2}|\nabla u_0|^2 + \frac{1}{2}|u_1|^2 + \frac{1}{p+1}|u_0|^{p+1}\right]dx < + \infty.
 \]
 Then we have
 \begin{itemize}
  \item [(a)] The corresponding solution $u$ to (CP1) satisfies a weighted Morawetz estimate
  \[
   \int_{-\infty}^\infty \int_{\Rm^d} (|t|+|x|)^\kappa \left(\frac{|\slashed{\nabla} u|^2}{|x|} + \lambda_d \frac{|u|^2}{|x|^3} + \frac{|u|^{p+1}}{|x|}\right) dx dt \lesssim_{d,p,\kappa} E_\kappa(u_0,u_1).
  \]
  In addition, the inward/outward energy satisfies the following decay estimates
 \begin{align*}
  &E_-(t) \lesssim_{d,p,\kappa} |t|^{-\kappa}, \; t>0;& &E_-(t) \in L^{1/\kappa}([0,\infty));&\\
  &E_+(t) \lesssim_{d,p,\kappa} |t|^{-\kappa}, \; t<0;& &E_+(t) \in L^{1/\kappa}((-\infty,0]).&
 \end{align*}
 It immediately follows that $u \in L^{(p+1)/\kappa} L^{p+1} (\Rm \times \Rm^d)$ since the definition of inward/outward energy implies that $\|u(\cdot,t)\|_{L^{p+1}(\Rm^d)}^{p+1} \lesssim E_\pm(t)$.
  \item[(b)] If $d,p,\kappa$ also satisfy $4\leq d\leq 8$, $p>p_c(d)$ and
  \[
   \kappa \geq \kappa_0(d,p) = \frac{(d+2)(d+3)-(d+3)(d-2)p}{(d-1)(d+3)-(d+1)(d-3)p} = \frac{p_e(d) - p}{(p_e(d)-p)+\frac{3(d-1)}{(d-2)(d+3)}(p-p_c(d))},
  \]
  then the corresponding solution $u$ to (CP1) scatters in both two time directions. More precisely, there exists $(v_0^\pm ,v_1^\pm) \in (\dot{H}^1\cap \dot{H}^{s_p}) \times (L^2\cap \dot{H}^{s_p-1})$, so that 
 \[
  \lim_{t \rightarrow \pm \infty} \left\|\begin{pmatrix} u(\cdot,t)\\ \partial_t u(\cdot,t)\end{pmatrix} - 
  \mathbf{S}_L (t)\begin{pmatrix}v_0^\pm \\ v_1^\pm\end{pmatrix}\right\|_{\dot{H}^s \times \dot{H}^{s-1}(\Rm^d)} = 0, \; \forall s\in [s_p,1].
 \]
 Here $\mathbf{S}_L (t)$ is the linear wave propagation operator. 
 \end{itemize}
\end{theorem}

\begin{remark} 
 Let us compare our scattering results in higher dimensions with those in 3-dimensional case. The minimal decay rate in 3-dimensional case is $\kappa_0(p)= (5-p)/2$, as given in \cite{shen3dnonradial}. This is consistent with the formula of $k_0(d,p)$ given above. The only difference between 3 and higher dimensional case is that the endpoint $\kappa = \kappa_0(d,p)$ is still allowed in the higher dimensional case. This is because in the proof we have to apply an interpolation between $L^{(p+1)/\kappa_0(d,p)} L^{p+1}$ and $L^2 L^\frac{2d}{d-3}$ spaces if we wish to include the endpoint case. This argument works perfectly in dimension 4 or higher. However, the pair $(2,\infty)$ is forbidden in the Strichartz estimates of 3D wave equation. We have to use $L^{2^+} L^{\infty^-}$ instead, which makes it necessary to assume $\kappa>\kappa_0(3,p)$. 
\end{remark}
\begin{remark}
 We can also prove a scattering result in the conformal case $p = p_c(d)$, $4\leq d \leq 9$. If the initial data $(u_0,u_1) \in \dot{H}^{1/2} \times \dot{H}^{-1/2}(\Rm^d)$ also satisfy
 \[
  \int_{\Rm^d} |x|\left(\frac{1}{2}|\nabla u_0(x)|^2 + \frac{1}{2}|u_1(x)|^2 + \frac{d-1}{2(d+1)}|u_0(x)|^{\frac{2(d+1)}{d-1}}\right) dx < +\infty, 
 \]
 then the corresponding solution $u$ to (CP1) must be a global-in-time solution and scatter in both time directions with $\|u\|_{W_d(\Rm)}\doteq\|u\|_{L^{\frac{2(d+1)}{d-1}} L^{\frac{2(d+1)}{d-1}}(\Rm \times\Rm^d)} < +\infty$. The proof is similar to the 3-dimensional case. Let us temporally assume the maximal lifespan of $u$ to be $(-T_-,T_+)$. On one hand, Remark \ref{gamma 1 remark} gives an integral estimate of $u$ when $|x|<2t$:
 \[
  \int_0^{T_+} \int_{\Rm^d} \frac{t}{|x|} |u(x,t)|^{\frac{2(d+1)}{d-1}} dx dt < +\infty \quad \Rightarrow \int_0^{T_+} \int_{|x|<2t} |u(x,t)|^{\frac{2(d+1)}{d-1}} dx dt < +\infty. 
 \]
Here we apply a smooth center cut-off technique when necessary. On the other hand, when $|x|>t+R$, we also have an inequality 
\[
 \int_{0}^{T_+} \int_{|x|>R+t} |u(x,t)|^{\frac{2(d+1)}{d-1}} dx dt < +\infty,
\]
as long as the radius $R$ is sufficiently large, by scattering of small solution, finite speed of propagation and a smooth cut-off technique. Combining these two integral estimates we are able to prove $\|u\|_{W_d([0,T_+))}<+\infty$. This immediately gives the global existence and scattering by scattering criterion Proposition \ref{scattering criterion}. Please refer to \cite{shen3dnonradial} for more details. 
\end{remark}

\section{Preliminary Results}

\subsection{Notations}

\paragraph{Derivatives} In this work the notation $\nabla u$ represents the gradient of $u$ with respect to the spatial variables 
\[
 \nabla u = \left(\frac{\partial u}{\partial x_1}, \frac{\partial u}{\partial x_2}, \cdots, \frac{\partial u}{\partial x_d}\right) \in \Rm^d.
\]
We define $\slashed{\nabla} u$ to be the covariant derivative of $u$ on the sphere centred at the origin with a fixed radius $|x|$. 
\[
 \slashed{\nabla} u = \nabla u - \left(\frac{x}{|x|}\cdot \nabla u\right)\frac{x}{|x|} \in \Rm^d.
\]
For convenience we also use the following notations for derivatives of $u$
\begin{align*}
 &u_r = \frac{\partial u}{\partial r};& & u_{i} = \frac{\partial u}{\partial x_i}; & &u_{ij} = \frac{\partial^2 u}{\partial x_i \partial x_j}.&
\end{align*}

\paragraph{Vectors} We use the notation $\vec{e}_i$ for the directional vector of $x_i$-axis and $\vec{e}$ for the directional vector of $t$-axis. We may use the notation $\vec{x}$ for the vector $x=(x_1,x_2,\cdots,x_d) \in \Rm^d$ when we emphasize that we are working with a vector. If necessary all the vectors in $\Rm_x^d$ are also treated as vectors in $\Rm_x^d \times \Rm_t$ by the natural embedding $x \rightarrow (x,0) \in \Rm_x^d \times \Rm_t$. For example, we may use the following notations
\begin{align*}
 &\nabla u = \left(\frac{\partial u}{\partial x_1}, \frac{\partial u}{\partial x_2}, \cdots, \frac{\partial u}{\partial x_d}, 0\right),& &\vec{x} =(x_1,x_2,\cdots,x_d,0);&
\end{align*}
if they appear in the expression of a vector field in $\Rm_x^d \times \Rm_t$.

\subsection{Technical Lemmata}

\begin{lemma} \label{identity of w u energy}
Let $u \in \dot{H}^{1}(\Rm^d)$ with $d \geq 3$. Then we have
\[
 \int_{\Rm^3} \left(\left|\mathbf{L} u\right|^2+ \lambda_d \cdot \frac{|u|^2}{|x|^2}\right) dx = \int_{\Rm^3} |u_r|^2 dx.
\]
\end{lemma}
\begin{proof}
We first consider the integral over annulus $\{x\in \Rm^d: a<|x|<b\}$: 
\begin{align}
 & \int_{a<|x|<b}  \left(\left|\mathbf{L} u\right|^2 + \lambda_d \frac{|u|^2}{|x|^2}\right) dx \nonumber\\
   =  &\int_{{\mathbb S}^{d-1}} \int_a^b \left(\left|u_r + \frac{d-1}{2}\cdot\frac{u}{r}\right|^2+ \lambda_d \frac{|u|^2}{r^2} \right) r^{d-1} \,dr d\Theta \nonumber\\
  = &\int_{{\mathbb S}^{d-1}} \int_a^b \left[r^{d-1} |u_r|^2 + \frac{d-1}{2}\partial_r (r^{d-2}|u|^2)\right] \,dr d\Theta \nonumber\\
    = &\int_{a<|x|<b} |u_r|^2 dx + \frac{d-1}{2b} \int_{|x|=b} |u|^2 d\sigma_b (x) - \frac{d-1}{2a} \int_{|x|=a} |u|^2 d\sigma_a (x). \label{relationship of w and u}
\end{align}
Here $\sigma_r$ represents regular measure on sphere of radius $r$. Next we use Hardy's inequality and obtain 
\[
 \int_0^\infty \left(\frac{1}{r^2} \int_{|x|=r} |u(x)|^2 d\sigma_r(x)\right) dr = \int_{\Rm^d} \frac{|u(x)|^2}{|x|^2} dx \lesssim\|u\|_{\dot{H}^1}^2 < +\infty.
\]
As a result we have
\begin{align}
 &\liminf_{r \rightarrow 0^+} \frac{1}{r} \int_{|x|=r} |u(x)|^2 d\sigma_r(x) = 0;& &\liminf_{r \rightarrow +\infty} \frac{1}{r} \int_{|x|=r} |u(x)|^2 d\sigma_r(x) = 0.& \label{limit of two sides}
\end{align}
Finally we may make $a \rightarrow 0^+$ and $b \rightarrow +\infty$ in identity \eqref{relationship of w and u} with these two limits in mind to finish the proof.
\end{proof}

\begin{remark} \label{relationship of u w energy}
 If we use one limit at a time in identity \eqref{relationship of w and u}, we obtain the following identities for any $\dot{H}^1(\Rm^d)$ function $u$ and any radius $R>0$
 \begin{align*}
  & \int_{|x|<R}  \left(\left|\mathbf{L} u(x)\right|^2+\lambda_d \frac{|u(x)|^2}{|x|^2}\right)  dx = \int_{|x|<R} |u_r|^2 dx + \frac{d-1}{2R} \int_{|x|=R} |u(x)|^2 d\sigma_R(x);\\
  & \int_{|x|>R}  \left(\left|\mathbf{L} u(x)\right|^2+\lambda_d \frac{|u(x)|^2}{|x|^2}\right)  dx = \int_{|x|>R} |u_r|^2 dx -  \frac{d-1}{2R} \int_{|x|=R} |u(x)|^2 d\sigma_R(x).
 \end{align*}
 This implies for any $\kappa > 0$ and $(u_0,u_1) \in \dot{H}^1 \times L^2$,
 \begin{align*}
   & \int_{\Rm^d} |x|^\kappa \left[\frac{1}{4}|\mathbf{L}_+ (u_0,u_1)|^2 + \frac{1}{4}|\mathbf{L}_- (u_0,u_1)|^2 + \frac{\lambda_d}{2} \frac{|u_0|^2}{|x|^2} + \frac{1}{2}|\slashed{\nabla} u_0|^2 +  \frac{1}{p+1}|u_0|^{p+1}\right] dx \\
   = &  \int_{\Rm^d} |x|^{\kappa} \left[\frac{|\slashed{\nabla} u_0|^2}{2} + \frac{|u_1|^2}{2} + \frac{|u_0|^{p+1}}{p+1}\right] dx
   + \int_0^\infty \frac{\kappa}{2} R^{\kappa-1} \left( \int_{|x|>R} \left(\left|\mathbf{L} u_0\right|^2 + \lambda_d \frac{|u_0|^2}{|x|^2}\right)dx\right) dR\\
   \leq & \int_{\Rm^d} |x|^{\kappa} \left[\frac{|\slashed{\nabla} u_0|^2}{2} + \frac{|u_1|^2}{2} + \frac{|u_0|^{p+1}}{p+1}\right] dx
   +  \int_0^\infty \frac{\kappa}{2} R^{\kappa-1} \left( \int_{|x|>R} |\partial_r u_0(x)|^2 dx\right) dR\\
   = & \int_{\Rm^d} |x|^\kappa \left[\frac{1}{2}|\nabla u_0|^2+ \frac{1}{2}|u_1|^2 + \frac{1}{p+1}|u_0|^{p+1}\right]  dx.
 \end{align*}
\end{remark}

\subsection{Strichartz estimates and local theory}

\paragraph{Strichartz estimates} We first recall the generalized Strichartz estimates of wave equation, which play a key role in the local theory. Please see Proposition 3.1 of Ginibre-Velo \cite{strichartz}. Here we use the Sobolev version. 

\begin{proposition}[Strichartz estimates]\label{Strichartz estimates} 
 Let $2\leq q_1,q_2 \leq \infty$, $2\leq r_1,r_2 < \infty$ and $\rho_1,\rho_2,s\in \Rm$ be constants with
 \begin{align*}
  &\frac{2}{q_i} + \frac{d-1}{r_i} \leq \frac{d-1}{2},& &(q_i,r_i)\neq \left(2,\frac{2(d-1)}{d-3}\right),& &i=1,2;& \\
  &\frac{1}{q_1} + \frac{d}{r_1} = \frac{d}{2} + \rho_1 - s;& &\frac{1}{q_2} + \frac{d}{r_2} = \frac{d-2}{2} + \rho_2 +s.&
 \end{align*}
 Assume that $u$ is the solution to the linear wave equation
\[
 \left\{\begin{array}{ll} \partial_t u - \Delta u = F(x,t), & (x,t) \in \Rm^d \times [0,T];\\
 u|_{t=0} = u_0 \in \dot{H}^s; & \\
 \partial_t u|_{t=0} = u_1 \in \dot{H}^{s-1}. &
 \end{array}\right.
\]
Then we have
\begin{align*}
 \left\|\left(u(\cdot,T), \partial_t u(\cdot,T)\right)\right\|_{\dot{H}^s \times \dot{H}^{s-1}} & +\|D_x^{\rho_1} u\|_{L^{q_1} L^{r_1}([0,T]\times \Rm^d)} \\
 & \leq C\left(\left\|(u_0,u_1)\right\|_{\dot{H}^s \times \dot{H}^{s-1}} + \left\|D_x^{-\rho_2} F(x,t) \right\|_{L^{\bar{q}_2} L^{\bar{r}_2} ([0,T]\times \Rm^d)}\right).
\end{align*}
Here the coefficients $\bar{q}_2$ and $\bar{r}_2$ satisfy $1/q_2 + 1/\bar{q}_2 = 1$, $1/r_2 + 1/\bar{r}_2 = 1$. The constant $C$ does not depend on $T$ or $u$. 
\end{proposition}

\begin{definition}
 We say $(q_1,r_1)$ is an $s$-admissible pair if the constants $q_1,r_1,s$ and $\rho_1 = 0$ satisfy the conditions given in Proposition \ref{Strichartz estimates}. For example, the pair $\left(\frac{(d+1)(p-1)}{2}, \frac{(d+1)(p-1)}{2}\right)$ is $s_p$-admissible as long as $p \in [p_c(d),p_e(d)]$. In particular, when $p=p_c(d)$, the pair $\left(\frac{2(d+1)}{d-1}, \frac{2(d+1)}{d-1}\right)$ is $1/2$-admissible. 
Because the space-time norms with these parameters will be frequently used in this work, we introduce the following notations for a time interval $I$
 \begin{align*}
  &\|u\|_{S_{d,p}(I)} = \|u\|_{L^{\frac{(d+1)(p-1)}{2}} L^{\frac{(d+1)(p-1)}{2}}(I \times \Rm^d)};&
  &\|u\|_{W_{d}(I)} = \|u\|_{L^{\frac{2(d+1)}{d-1}} L^{\frac{2(d+1)}{d-1}}(I \times \Rm^d)};& 
 \end{align*}
\end{definition}

\paragraph{Chain rule} We also need the following ``chain rule'' for fractional derivatives. Please refer to Christ-Weinstein \cite{fchain}, Kato \cite{kato}, Kenig et al. \cite{fchain2}, Staffilani \cite{fchain3} and Taylor \cite{fchain4} for more details. 
\begin{lemma} \label{chain rule}
 Assume a $C^1$ function $F$ satisfies $F(0)=0$ and
 \begin{align*}
  |F'(\theta u +(1-\theta)v)| \leq C(|F'(u)|+|F'(v)|), \quad \forall \theta \in [0,1].
 \end{align*}
Then we have
 \[
  \|D^\alpha F(u)\|_{L^q (\Rm^d)} \lesssim \|D^\alpha u\|_{L^{q_1} (\Rm^d)} \|F'(u)\|_{L^{q_2} (\Rm^d)}
 \]
 for $0<\alpha<1$ and $1/q = 1/q_1+1/q_2$, $1<q,q_1,q_2<\infty$.
\end{lemma}

\paragraph{Local theory} We may combine the Strichartz estimates with a fixed point argument to develop a local theory of (CP1), with initial data in either the critical Sobolev space $\dot{H}^{s_p} \times \dot{H}^{s-p-1}$ or the energy space $\dot{H}^1 \times L^2$. Please see Kapitanski \cite{loc1} and Lindblad-Sogge \cite{ls}, for instance, for more results and details about the local theory.

\begin{proposition}[Existence and uniqueness] \label{local existence}
Assume $d\geq 3$ and $p\in [p_c(d),p_e(d)]$. For any initial data $(u_0,u_1) \in \dot{H}^{s_p} \times \dot{H}^{s_p-1}(\Rm^d)$, there exists a unique solution $u$ to (CP1) with a maximal lifespan $(-T_-,T_+)$ so that
\begin{itemize}
  \item $(u(\cdot,t), \partial_t u(\cdot,t)) \in C((-T_-,T_+); \dot{H}^{s_p}\times \dot{H}^{s_p-1})$;
  \item The norms $\|u\|_{S_{d,p}([a,b])}$ and $\|D_x^{s_p-1/2} u\|_{L^{\frac{2(d+1)}{d-1}} L^{\frac{2(d+1)}{d-1}} ([a,b]\times \Rm^3)}$ are finite for any time interval $[a,b]$ with $-T_-<a<b<T_+$. 
 \end{itemize}
\end{proposition} 
\begin{proposition}[Scattering criterion] \label{scattering criterion}
Assume $d\geq 3$ and $p\in [p_c(d),p_e(d)]$. If a solution $u$ to (CP1) with a maximal lifespan $(-T_-,T_+)$ satisfies $\|u\|_{S_{d,p}([0,T_+))}<+\infty$,
then $T_+=+\infty$ and the solution $u$ scatters in the positive time direction.
\end{proposition}

\begin{lemma} \label{lemma local in energy space}
Assume that $(d,p)$ satisfies assumption (A1). Let $(u_0,u_1) \in \dot{H}^1 \times L^2(\Rm^d)$ be initial data. Then the Cauchy problem (CP1) has a unique solution $u$ in the time interval $[0,T]$ with $(u(\cdot,t),u_t(\cdot,t)) \in C([0,T];\dot{H}^1 \times L^2(\Rm^d))$ and $u \in L^{\frac{2p}{(d-2)p-d}}L^{2p} ([0,T]\times \Rm^d)$. Here the minimal time length of existence 
\[
  T = C(d,p) \|(u_0,u_1)\|_{\dot{H}^1\times L^2}^{\frac{-2(p-1)}{(d+2)-(d-2)p}}.
\]
\end{lemma}

\paragraph{Global existence} Let $u$ be a solution to (CP1) with a finite energy $E$. If $u$ is defined at time $t_0>0$, then it can also be defined for all time in the  interval $[t_0,t_0+T]$, with 
\[
 T = C(d,p) \|(u(\cdot,t_0), u_t(\cdot,t_0))\|_{\dot{H}^1 \times L^2}^{\frac{-2(p-1)}{(d+2)-(d-2)p}} \geq C(d,p) (2E)^{-\frac{p-1}{(d+2)-(d-2)p}},
\]
according to Lemma \ref{lemma local in energy space}. Since the lower bound of $T$ does not depend on $t_0$, the solution $u$ can be defined for all time $t\in \Rm^+$. The same argument also works in the negative time direction since the wave equation is time-reversible. 

\begin{proposition} \label{global existence finite energy}
Assume that $(d,p)$ satisfies assumption (A1). If $u$ is a solution to (CP1) with a finite energy, then $u$ is defined for all time $t \in \Rm$. 
\end{proposition}

\begin{remark} \label{local L1L2}
The assumption (A1) is actually a combination of the following two inequalities 
\begin{align*}
 &p_c(d) \leq p < p_e(d);\quad (a1)& &p \leq 1+\frac{3}{d-3}.\quad (a2)&
\end{align*}
The author would like to explain why we make the technical assumption (a2). In fact, assumption (a2) is equivalent to $\frac{2p}{(d-2)p-d}\geq 2$, which is a necessary condition for $(\frac{2p}{(d-2)p-d}, 2p)$ to be an admissible pair, thus (a2) is essential to Lemma \ref{lemma local in energy space}.  If $(d,p)$ satisfies (A1) and $u$ is a solution to (CP1) with a finite energy, then we may apply Lemma \ref{lemma local in energy space} to obtain 
 \[
  u \in L_{loc}^{\frac{2p}{(d-2)p-d}} L^{2p} (\Rm\times \Rm^d) \hookrightarrow L_{loc}^{p} L^{2p} (\Rm\times \Rm^d) \quad \Rightarrow \quad F(u) \in L_{loc}^1 L^2(\Rm \times \Rm^d).
 \] 
Please note that the inequality $\frac{2p}{(d-2)p-d} \geq p$ always holds as long as $p \leq p_e(d)$. The fact $u \in L_{loc}^p L^{2p}$ plays an important role in smooth approximation argument, as shown in Remark \ref{smooth approximation Morawetz} and \ref{smooth approximation flux formula}. If $p > 1+\frac{3}{d-3}$, however, we can no longer follow a similar argument to show $u \in L_{loc}^p L^{2p}$ because we can never find $s \leq 1$ and $q \geq 2$ so that $(q,2p)$ is an $s$-admissible pair in this case. 
\end{remark}

\subsection{Morawetz estimates} \label{sec: Morawetz subsection}
The topic of this subsection is a Morawetz inequality of wave equation. This immediately gives a few global integral estimates, which are crucial to the development of our inward/outward energy theory. Our Morawetz inequality is similar to the Morawetz inequality given by Perthame and Vega in the final section of their work \cite{benoit}. The difference is that our new version of Morawetz inequality comes with two additional terms, i.e. the integrals of $|\slashed{\nabla} u|^2 /|x|$ and $|u|^2/|x|^3$, in the left hand. 
\begin{theorem} 
Assume that $(d,p)$ satisfies (A1). Let $u$ be a solution to (CP1) with a finite energy $E$. Then we have the following inequality for any $R>0$. Here $\sigma_R$ is the regular surface measure of the sphere $|x|=R$. 
\begin{align}
 & \frac{1}{2R}\int_{-\infty}^\infty \!\int_{|x|<R}\left(|\nabla u|^2+|u_t|^2+\frac{(d\!-\!1)(p\!-\!1)\!-\!2}{p+1}|u|^{p+1}\right) dx dt + \frac{d-1}{4R^2} \int_{-\infty}^\infty \!\int_{|x|=R} |u|^2 d\sigma_R dt  \nonumber \\
 & \quad + \int_{-\infty}^\infty\! \int_{|x|>R} \left(\frac{|\slashed{\nabla} u|^2}{|x|}+ \lambda_d \frac{|u|^2}{|x|^3}+ \frac{(d-1)(p-1)}{2(p+1)} \frac{|u|^{p+1}}{|x|}\right) dx dt  \leq 2E. \label{morawetz1}
\end{align}
\end{theorem}
\begin{proof}
Perthame and Vega sketch a proof of their version of Morawetz inequality in the final section of their work \cite{benoit}. We will follow almost the same argument but with two improvements, which enable us to insert two more terms mentioned above to the left hand of the inequality. We will also give more details for reader's convenience. Given a positive constant $R$, we define two radial functions $\Psi$ and $\varphi$ by
 \begin{align*}
  &\nabla \Psi = \left\{\begin{array}{ll} x, & \hbox{if} \; |x|\leq R; \\ Rx/|x|, & \hbox{if}\; |x|\geq R;\end{array}\right.&
  &\varphi = \left\{\begin{array}{ll} 1/2, & \hbox{if} \; |x|\leq R; \\ 0, & \hbox{if}\; |x|> R.\end{array}\right.&
 \end{align*}
 Since $u$ is defined for all time $t \in \Rm$, we may also define a function on $\Rm$
 \[
  \mathcal{E}(t) = \int_{\Rm^d}  u_t(x,t) \left(\nabla u(x,t)\cdot \nabla \Psi + u(x,t)\left(\frac{\Delta \Psi}{2} - \varphi\right)\right) dx.
 \]
 We may differentiate $\mathcal{E}$, utilize the equation $u_{tt} - \Delta u =  -|u|^{p-1}u$, apply integration by parts and obtain
 \begin{align*}
   -\mathcal{E}'(t) & = \int_{\Rm^d} \left(\sum_{i,j=1}^d u_i \Psi_{ij} u_j-\varphi |\nabla u|^2 + \varphi |u_t|^2 \right) dx + \frac{1}{4} \int_{\Rm^d} \nabla (|u|^2)\cdot  \nabla \left(\Delta \Psi - 2\varphi\right) dx\\
   & \qquad \qquad + \int_{\Rm^d} |u|^{p+1}\left(\frac{p-1}{2(p+1)}\Delta \Psi - \varphi\right) dx = I_1 + I_2 + I_3.
 \end{align*}
 Here we have
 \begin{align*}
  &\Psi_{ij} = \left\{\begin{array}{ll} \delta_{ij}, & \hbox{if} \; |x|< R; \\ \frac{R\delta_{ij}}{|x|} - \frac{R x_i x_j}{|x|^3}, & \hbox{if}\; |x|> R;\end{array}\right.&
  &\Delta \Psi = \left\{\begin{array}{ll} d, & \hbox{if} \; |x|< R; \\ R(d-1)/|x|, & \hbox{if}\; |x|> R;\end{array}\right.&
 \end{align*}
 \begin{align*}
  \Delta \Psi - 2\varphi = \left\{\begin{array}{ll} d-1, & \hbox{if} \; |x|\leq R; \\ R(d-1)/|x|, & \hbox{if}\; |x|\geq R;\end{array}\right. \in C(\Rm^d).
 \end{align*}
when $|x|>R$, we may calculate
\begin{align*}
 \sum_{i,j=1}^d u_i \Psi_{ij} u_j = \sum_{i,j=1}^d u_i \left(\frac{R\delta_{ij}}{|x|} - \frac{R x_i x_j}{|x|^3}\right) u_j
= \frac{R}{|x|}|\nabla u|^2 - \frac{R|\nabla u \cdot x|^2}{|x|^3}
= \frac{R}{|x|}|\slashed{\nabla} u|^2.
\end{align*} 
Thus we have 
\begin{align}
 I_1 =  \frac{1}{2} \int_{|x|<R} \!\!\left(|\nabla u|^2 + |u_t|^2\right) dx + R\int_{|x|>R} \frac{|\slashed{\nabla} u|^2}{|x|} dx. \label{Morawetz contribution 1}
\end{align}
The last term in the equality above, which vanishes for radial solutions, is discarded in the original argument. A basic computation shows
\begin{equation}
 I_3 = \frac{(d\!-\!1)(p\!-\!1)\!-\!2}{2(p+1)}\int_{|x|<R} |u|^{p+1} dx  + \frac{(d-1)(p-1)R}{2(p+1)} \int_{|x|>R} \frac{|u|^{p+1}}{|x|} dx. \label{Morawetz contribution 2}
\end{equation}
Our second improvement is about the estimate of $I_2$. The original paper uses the following inequality
\[
  \frac{1}{4} \int_{\Rm^d} \nabla (|u|^2) \nabla \left(\Delta \Psi - 2\varphi\right) dx \geq \frac{d-1}{4R} \int_{|x|=R} |u|^2 d\sigma_R (x).
\]
It turns out that we may improve this inequality in higher dimensions $d \geq 4$. Let us calculate the left hand carefully. 
\begin{align}
 I_2 & = \frac{1}{4}\int_{\Rm^d} \nabla (|u|^2) \cdot \nabla \left(\Delta \Psi - 2\varphi\right) dx \nonumber\\
 & = \frac{1}{4}\int_{|x|>R} \nabla (|u|^2) \cdot \frac{-R (d-1)x}{|x|^3} dx \nonumber\\
 & = \frac{1}{4} \int_{|x|>R} \left[\hbox{div}\left(|u|^2 \cdot \frac{-R(d-1) x}{|x|^3}\right) + 4\lambda_d \frac{R}{|x|^3} |u|^2\right]dx \nonumber\\
 & = \frac{d-1}{4R} \int_{|x|=R} |u|^2 d\sigma_R(x) + \lambda_d R \int_{|x|>R} \frac{|u|^2}{|x|^3} dx. \label{Morawetz contribution 3}
\end{align}
Since $-\mathcal{E}'(t) = I_1 + I_2 + I_3$, we have 
\begin{equation}
 \int_{t_1}^{t_2} (I_1+I_2+I_3) dt = \mathcal{E}(t_1) - \mathcal{E}(t_2). \label{Morawetz integral identity}
\end{equation}
A uniform upper bound of $|\mathcal{E}(t)|$ can be found by
\begin{align*}
 R |\mathcal{E}(t)| &\leq \frac{1}{2}\int_{\Rm^d} \left( R^2|u_t(x,t)|^2 +  \left|\nabla u(x,t)\cdot \nabla \Psi + u(x,t)\left(\frac{\Delta \Psi}{2} - \varphi\right)\right|^2 \right) dx\\
 & = \frac{1}{2}\int_{\Rm^d} \left(R^2 |u_t|^2 + |\nabla u\cdot \nabla \Psi|^2 + \left(\frac{\Delta \Psi}{2} - \varphi\right) \nabla (|u|^2) \cdot \nabla \Psi + \left(\frac{\Delta \Psi}{2} - \varphi\right)^2 |u|^2 \right) dx\\
 & \leq \frac{1}{2} \int_{\Rm^d} \left[R^2 |u_t|^2 + R^2 |\nabla u|^2 - \hbox{div} \left(\left(\frac{\Delta \Psi}{2} - \varphi\right)\nabla \Psi \right) |u|^2 + \left(\frac{\Delta \Psi}{2} - \varphi\right)^2 |u|^2\right] dx
\end{align*}
A basic calculation shows 
\[
 \hbox{div} \left(\left(\frac{\Delta \Psi}{2} - \varphi\right)\nabla \Psi \right) = 
 \left\{\begin{array}{ll} d(d-1)/2, & \hbox{if} \; |x|\leq R; \\ \frac{(d-1)(d-2)R^2}{2|x|^2}, & \hbox{if}\; |x|\geq R;\end{array}\right\} \geq \left(\frac{\Delta \Psi}{2} - \varphi\right)^2.
\]
Thus we have 
\[
 R |\mathcal{E}(t)| \leq R^2 E \quad \Rightarrow \quad |\mathcal{E}(t)| \leq RE.
\]
Plugging this upper bound into the integral identity \eqref{Morawetz integral identity}, we obtain 
\begin{align} \label{finite time estimate M}
 \int_{t_1}^{t_2} (I_1+I_2+I_3) dt \leq 2RE  
\end{align}
for all $t_1<t_2$. Next we make $t_1\rightarrow -\infty$, $t_2\rightarrow +\infty$, divide both sides by $R$ and obtain an inequality
\begin{equation*} 
 \frac{1}{R} \int_{-\infty}^\infty (I_1+I_2+I_3) dt \leq 2E.
\end{equation*}
Finally we plug in the expressions of $I_1$, $I_2$, $I_3$ given above into this inequality and verify that this is exactly inequality \eqref{morawetz1} we want to prove. 
\end{proof}

\begin{remark} \label{smooth approximation Morawetz}
 In the proof above we calculate as thought $u \in C^2(\Rm^d\times \Rm)$ decays sufficiently fast when $|x| \rightarrow +\infty$. In order to deal with the general case, we may use the following smooth approximation and cut-off techniques. Let $\phi_d$ be a smooth, radial, nonnegative, compactly supported function defined on $\Rm^d$ with $\|\phi_d\|_{L^1(\Rm^d)} = 1$. We define 
  \begin{align*}
   &\varphi(x,t) = \phi_d(x) \phi_1(t), \; (x,t)\in \Rm^d \times \Rm;& &\varphi_\eps (x,t) = \eps^{-d-1} \varphi(x/\eps,t/\eps), \; \eps>0;&
   \end{align*}
  If $u$ is a solution to (CP1) with compactly supported initial data, then the smooth function $u_\eps = \phi_\eps \ast u$ solves the equation 
  \begin{align*}
   &\partial_t^2 u_\eps - \Delta u_\eps = F_\eps(x,t); & &F_\eps = \phi_\eps \ast (-|u|^{p-1}u).&
  \end{align*}
  The smooth functions $u_\eps$ and $F_\eps$ satisfy the energy estimates
  \begin{align*}
   & E_\eps(t) \doteq \int_{\Rm^d} \left(\frac{1}{2}|\nabla u_\eps(x,t)| + \frac{1}{2}|\partial_t u_\eps (x,t)|^2 + \frac{1}{p+1}|u_\eps|^{p+1}\right) dx \leq E;\\
   & \lim_{\eps \rightarrow 0^+} \|(u_\eps,\partial_t u_\eps)-(u,u_t)\|_{C([t_1,t_2]; \dot{H}^{1}\times L^2(\Rm^d))} =  0, \quad -\infty <t_1<t_2<+\infty. 
  \end{align*}
  We also have the following convergence of space-time norms
  \begin{align*}
   \lim_{\eps \rightarrow 0^+} & \left(\|F_\eps +|u|^{p-1} u\|_{L^1 L^2([t_1,t_2]\times \Rm^d)} + \|u_\eps - u\|_{L^p L^{2p}([t_1,t_2]\times \Rm^d)}\right)= 0;\\
   & \qquad \qquad \Rightarrow \lim_{\eps \rightarrow 0^+} \|F_\eps +|u_\eps|^{p-1} u_\eps\|_{L^1 L^2([t_1,t_2]\times \Rm^d)} = 0.
  \end{align*}
  Now we may substitute $u$ by $u_\eps$ and proceed as in the proof given above. Since $u_\eps$ is no longer a solution to (CP1), we need to insert an additional error term to the expression of $\mathcal{E}'_{\eps} (t)$ 
  \[
    Err = \int_{\Rm^d} (F_\eps + |u_\eps|^{p-1} u_\eps) \left(\nabla u_\eps(x,t)\cdot \nabla \Psi + u_\eps(x,t)\left(\frac{\Delta \Psi}{2} - \varphi\right)\right) dx,
  \]
  whose upper bound is given by
  \begin{align*}
   |Err| & \leq \|F_\eps + |u_\eps|^{p-1} u_\eps\|_{L^2(\Rm^d)} \left\|\nabla u_\eps(x,t)\cdot \nabla \Psi + u_\eps(x,t)\left(\frac{\Delta \Psi}{2} - \varphi\right)\right\|_{L^2(\Rm^d)}\\
   & \lesssim_d \|F_\eps + |u_\eps|^{p-1} u_\eps\|_{L^2(\Rm^d)} \left(R\|\nabla u_\eps\|_{L^2(\Rm^d)} + \left\|\frac{R}{|x|}u_\eps\right\|_{L^2(\Rm^d)}\right)\\
   & \lesssim_d R E^{1/2} \|F_\eps + |u_\eps|^{p-1} u_\eps\|_{L^2(\Rm^d)}.
  \end{align*}
  Thus the space-time integral estimate \eqref{finite time estimate M} can be rewritten as 
  \[
    \int_{t_1}^{t_2} (I_{1,\eps}+I_{2,\eps}+I_{3,\eps}) dt \leq RE_\eps(t_1) + RE_\eps(t_2) + CR E^{1/2} \|F_\eps + |u_\eps|^{p-1} u_\eps\|_{L^1 L^2([t_1,t_2] \times \Rm^d)}.
  \]
  Please note that the $L^1 L^2$ norm above vanishes as $\eps \rightarrow 0^+$. Thus we may make $\eps \rightarrow 0^+$, then let $t_1\rightarrow -\infty$, $t_2\rightarrow +\infty$ and finish the proof for compactly supported initial data. In the case of general initial data we may apply the standard cut-off techniques and utilize the finite speed of propagation of wave equation.
\end{remark} 

\begin{corollary}
 Let $u$ be a solution to (CP1) with a finite energy $E$. Then we have
 \begin{align*}
  \int_{-\infty}^\infty \int_{\Rm^d} \left(\frac{|\slashed{\nabla} u|^2}{|x|}+ \lambda_d \frac{|u|^2}{|x|^3}+ \frac{|u|^{p+1}}{|x|}\right) dx dt & \lesssim_{d,p} E; \\
  \int_{-\infty}^\infty \!\int_{|x|<R}\left(|\nabla u|^2+|u_t|^2+|u|^{p+1}\right) dx dt & \lesssim_{d,p} RE;\\
  \int_{-\infty}^\infty \!\int_{|x|=R} |u|^2 d\sigma_R(x) dt \lesssim_{p,d} R^2E.
 \end{align*}
\end{corollary}
\section{Inward/outward Energy Theory}

\subsection{Statement of energy flux formula} \label{statement eff}

\begin{figure}[h]
 \centering
 \includegraphics[scale=1.0]{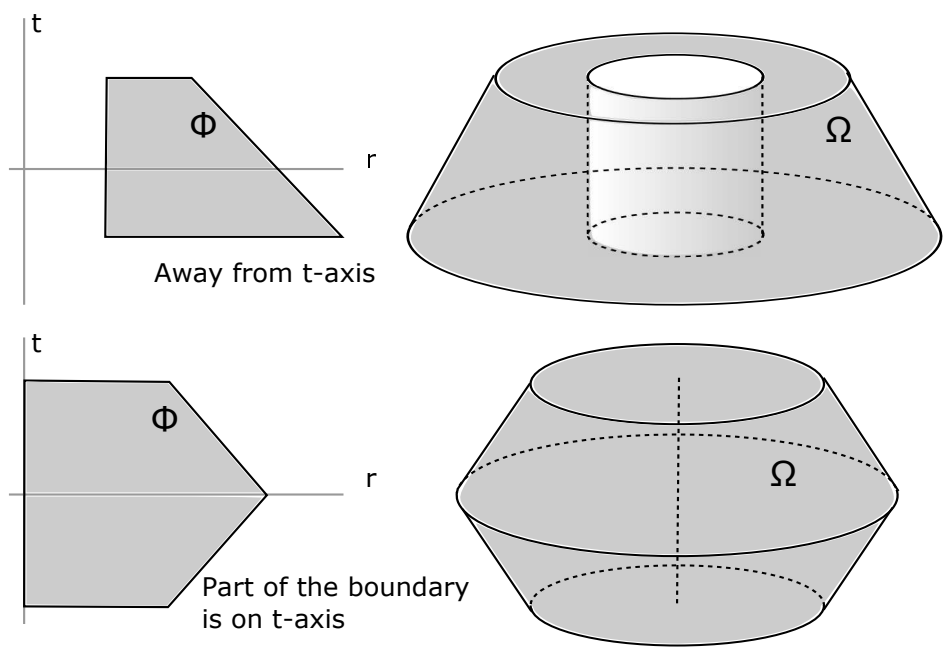}
 \caption{Two examples of regions} \label{figure regions}
\end{figure}

\paragraph{Regions in this work} For convenience we focus on radially symmetric regions as described below. Let $\Omega \subset \Rm^d \times \Rm$ be a region so that it can be expressed by
\[
 \Omega = \{(r\Theta, t)\in \Rm^d \times \Rm: (r,t)\in \Phi, \Theta \in \mathbb{S}^{d-1}\},
\]
if we use spherical coordinates $(r,\Theta)$ in $\Rm^d$. Here $\Phi \subset \Rm_r^+ \times \Rm_t$ is a bounded, closed region, whose boundary is a simple curve consisting of finite line segments paralleled to either $t$-axis, $r$-axis or $t \pm r=0$. Therefore the boundary surface $\partial \Omega$ consists of finite pieces of annulus, circular cylinders or cones. In addition we also allow a line segment of $t$-axis to be part of the boundary $\partial \Phi$. In this case the boundary surface $\partial \Omega$ contains a degenerate part, i.e. the same line segment of $t$-axis as mentioned above. Figure \ref{figure regions} shows two examples of these regions. One may also consider regions that are not radially symmetric. Please see Remark \ref{nonradial regions} below. 

\begin{proposition}[General Energy Flux] \label{energy flux formula}
 Assume that $(d,p)$ satisfies (A1). For connivence we first introduce the notation 
 \[
  e'(x,t) = \frac{1}{2}|\slashed{\nabla} u(x,t)|^2 + \frac{\lambda_d}{2}\cdot \frac{|u(x,t)|^2}{|x|^2} + \frac{1}{p+1}|u(x,t)|^{p+1} 
 \]
for the non-directional part of the integrand in the definition of energy then define two vector fields in $\Rm^d \times \Rm$:
 \begin{align*}
  \mathbf{V}_- = & \left[ -\frac{1}{4} \left|\mathbf{L}_+ u\right|^2 + \frac{1}{2} e'(x,t)\right]\frac{\vec{x}}{|x|} + \left[\frac{1}{4}\left|\mathbf{L}_+ u\right|^2 + \frac{1}{2}e'(x,t)\right]\vec{e};\\
   \mathbf{V}_+ = & \left[ +\frac{1}{4} \left|\mathbf{L}_- u\right|^2 - \frac{1}{2} e'(x,t)\right]\frac{\vec{x}}{|x|} + \left[\frac{1}{4}\left|\mathbf{L}_- u\right|^2 + \frac{1}{2}e'(x,t)\right]\vec{e};
 \end{align*}
Let $\Omega = \{(r\Theta, t)\in \Rm^d \times \Rm: (r,\Theta) \in \Phi, \Theta \in \mathbb{S}^{d-1}\}$ be a region as described above. If $\Phi$ also satisfies $\Phi \subset \Rm_r^+ \times \Rm_t$, then we have
 \[
  \int_{\partial \Omega} \mathbf{V}_{\pm} \cdot d\mathbf{S} = \pm  \iint_\Omega M(x,t) dx dt.
 \]
 Here the right hand is the Morawetz integral, whose integrand $M$ is the Morawetz density function 
 \[
  M(x,t) = \frac{1}{2} \cdot \frac{|\slashed{\nabla} u(x,t)|^2}{|x|} + \frac{\lambda_d}{2} \cdot \frac{|u(x,t)|^2}{|x|^3} + \frac{(d-1)(p-1)}{4(p+1)}\cdot \frac{|u(x,t)|^{p+1}}{|x|}.
 \]
 In addition, there exist a nonnegative, finite and continuous\footnote{Continuity means that the function $\mu((-\infty,x])$ is a continuous function of $x \in \Rm$.} Borel measure $\mu$ on $\Rm$ with $\mu(\Rm) \lesssim E$, which is determined by $u$ and independent to $\Omega$, so that if the line segment $(t_1,t_2)$ on the $t$-axis is a part of $\partial \Phi$, then the identity above still holds if we add $\mp c_d \cdot \mu([t_1,t_2])$ to the left hand side. Here $c_d$ is $(d-1)^2/16$ times the area of the unit sphere $\mathbb{S}^{d-1}$.
\end{proposition}

\paragraph{Surface integrals} For the reader's convenience we write down the details of surface integrals in the energy flux formula over different kinds of boundary surfaces $\Sigma$, as shown in table \ref{surface integrals in energy flux}. The physical interpretation of these surface integrals are similar to the 3-dimensional case. In general, we have 
\begin{itemize} 
 \item We can split the integral of $e'(x,t)$ into two parts: the integrals of $|\slashed{\nabla} u|^2$ and $|u^2|/|x|^2$ represent energy gain or loss by a combination of boundary effect and linear propagation; the integrals of $|u|^{p+1}$ represent the energy gain or loss by a combination of boundary and nonlinear effect. 
 \item The integrals of $|\mathbf{L}_\pm u|^2$ give the amount of energy moving through the surface $\Sigma$ due to linear propagation. 
 \item The measure $c_d \mu([t_1,t_2])$ is equal to the amount of energy carried by inward waves that move through the origin and change to outward waves during the time period $[t_1,t_2]$.
\end{itemize} 
Furthermore, the Morawetz integral $\iint_\Omega M(x,t) dx dt$ gives the amount of inward energy transformed to outward energy in the given space-time region $\Omega$. More details about the physical interpretation and the way to apply energy flux formula can be found in \cite{shen3dnonradial}.

\begin{table}[h]
\caption{Surface integrals in energy flux formula}
\begin{center}
\begin{tabular}{|c|c|c|}\hline
 Boundary type & Inward Energy Case & Outward Energy Case\\
 \hline
 Horizontally $\uparrow$ & $\int_\Sigma \left(\frac{\left|\mathbf{L}_+ u\right|^2}{4} \!+\! \frac{e'(x,t)}{2}\right) dS$ & 
 $\int_\Sigma \left(\frac{\left|\mathbf{L}_- u\right|^2}{4} \!+\! \frac{e'(x,t)}{2}\right) dS$ \\
 \hline
 Horizontally $\downarrow$ & $-\int_\Sigma \left(\frac{\left|\mathbf{L}_+ u\right|^2}{4} \!+\! \frac{e'(x,t)}{2}\right) dS$ & 
 $-\int_\Sigma \left(\frac{\left|\mathbf{L}_- u\right|^2}{4} \!+\! \frac{e'(x,t)}{2}\right) dS$\\
 \hline
 $\!|x|\!=\!r_0$, Outward & $\int_\Sigma \left(-\frac{\left|\mathbf{L}_+ u\right|^2}{4} \!+\! \frac{e'(x,t)}{2}\right) dS$ & 
 $\int_\Sigma \left(\frac{\left|\mathbf{L}_- u\right|^2}{4} \!-\! \frac{e'(x,t)}{2}\right) dS$ \\
 \hline
$\!|x|\!=\!r_0$, Inward & $\int_\Sigma \left(\frac{\left|\mathbf{L}_+ u\right|^2}{4} \!-\! \frac{e'(x,t)}{2}\right) dS$ & 
 $\int_\Sigma \left(-\frac{\left|\mathbf{L}_- u\right|^2}{4} \!+\! \frac{e'(x,t)}{2}\right) dS$ \\
 \hline
 $|x|=0$ & $-c_d \mu([t_1,t_2])$ & $c_d \mu([t_1,t_2])$ \\
 \hline
 Backward Cone$\uparrow$ & $\frac{1}{\sqrt{2}} \int_\Sigma e'(x,t) dS$ & $\frac{1}{2\sqrt{2}}\int_\Sigma |\mathbf{L}_- u|^2 dS$\\
 \hline
 Backward Cone$\downarrow$ & $-\frac{1}{\sqrt{2}} \int_\Sigma e'(x,t) dS$ & $-\frac{1}{2\sqrt{2}}\int_\Sigma |\mathbf{L}_- u|^2 dS$\\
 \hline
 Light Cone $\uparrow$ & $\frac{1}{2\sqrt{2}}\int_\Sigma |\mathbf{L}_+ u|^2 dS$ & $\frac{1}{\sqrt{2}} \int_\Sigma e'(x,t) dS$\\
 \hline
 Light Cone $\downarrow$ & $-\frac{1}{2\sqrt{2}}\int_\Sigma |\mathbf{L}_+ u|^2 dS$ & $-\frac{1}{\sqrt{2}} \int_\Sigma e'(x,t) dS$\\
  \hline
\end{tabular}
\end{center}
\label{surface integrals in energy flux}
\end{table}

\paragraph{Fluxes through light cones} For connivence we introduce the following notations.
\begin{definition} \label{def fluxes light cone}
 Let $C^-(s)$ and $C^+(\tau)$ be backward and forward light cones, respectively
 \begin{align*}
  & C^-(s) = \{(x,t)\in \Rm^d \times \Rm, |x|+t = s\};& &C^+(\tau) = \{(x,t)\in \Rm^d \times \Rm: t-|x| = \tau\}.&
 \end{align*}
 We then define $Q$'s to be the energy fluxes through these light cones. The upper index $\pm$ indicates whether the light cone is a forward one ($+$) or a backward one ($-$); while the lower index $\pm$ indicates whether this is the flux of inward energy ($-$) or outward energy ($+$).
\begin{align*} 
 Q_-^- (s) & = \frac{1}{\sqrt{2}} \int_{C^-(s)} \left(\frac{1}{2}\left|\slashed{\nabla} u\right|^2 + \frac{\lambda_d}{2} \cdot \frac{|u|^2}{|x|^2} + \frac{1}{p+1}|u|^{p+1} \right) dS;\\
 Q_+^-(s) & = \frac{1}{2\sqrt{2}} \int_{C^-(s)} |\mathbf{L}_- u|^2 dS;\\
 Q_-^+(\tau) & = \frac{1}{2\sqrt{2}} \int_{C^+(\tau)} |\mathbf{L}_+ u|^2 dS;\\
 Q_+^+(\tau) & = \frac{1}{\sqrt{2}} \int_{C^+(\tau)} \left(\frac{1}{2}\left|\slashed{\nabla} u\right|^2 + \frac{\lambda_d}{2} \cdot \frac{|u|^2}{|x|^2} + \frac{1}{p+1}|u|^{p+1} \right) dS.
\end{align*}
\end{definition}

\begin{remark} All the energy fluxes $Q$'s defined above is dominated by the energy $E$. In fact, the sum $Q_-^+(\tau) + Q_+^+(\tau)$ is exactly the flux of the full energy through the light cone $C^+(\tau)$, as shown below.
 \begin{align*}
  Q_-^+(\tau)+Q_+^+(\tau) & = \frac{1}{\sqrt{2}} \int_{C^+(\tau)} \left(\frac{1}{2}\left|\mathbf{L}_+ u\right|^2 + \frac{\lambda_d}{2}\frac{|u|^{2}}{|x|^2} +  \frac{1}{2}\left|\slashed{\nabla} u\right|^2  + \frac{1}{p+1}|u|^{p+1}\right) dS \\
  & = \int_{\Rm^d} \left(\frac{1}{2}\left|\mathbf{L} \tilde{u}\right|^2 + \frac{\lambda_d}{2}\frac{|\tilde{u}|^{2}}{|x|^2} +  \frac{1}{2}\left|\slashed{\nabla} \tilde{u}\right|^2  + \frac{1}{p+1}|\tilde{u}|^{p+1}\right) dx \\
  & = \int_{\Rm^d} \left(\frac{1}{2}\left|\partial_r \tilde{u}\right|^2 +  \frac{1}{2}\left|\slashed{\nabla} \tilde{u}\right|^2  + \frac{1}{p+1}|\tilde{u}|^{p+1}\right) dx \\
  & = \frac{1}{\sqrt{2}} \int_{C^+(\tau)} \left(\frac{1}{2}\left|(\partial_r + \partial_t)u\right|^2 +  \frac{1}{2}\left|\slashed{\nabla} u\right|^2  + \frac{1}{p+1}|u|^{p+1}\right) dS \leq E.
 \end{align*}
Here we use the notation $\tilde{u}(x) = u(x,\tau+|x|) \in \dot{H}^1(\Rm^d)$ and apply Lemma \ref{identity of w u energy}. The case $Q_-^-(s)+Q_+^-(s)$ is similar. 
\end{remark}

\subsection{Proof of Energy Flux Formula}

Now let us prove the energy flux formula. The proof consists of two major steps. First we apply Gauss' Formula to show that the energy flux formula holds for regions $\Omega$ away from the $t$-axis. Then we deal with other regions $\Omega$ by applying the formula on the region $\Omega_r = \Omega\cap \{(x,t)\in \Rm^d \times \Rm: |x|\geq r\}$ with an arbitrarily small positive constant $r$ and taking a limit process $r\rightarrow 0^+$. The second step does not depend on the dimension $d$, i.e. we may apply the same argument as in the case of dimension 3. Thus we omit the details of this argument in the work. Readers may look at \cite{shen3dnonradial} for details. The author would like to mention that if the solution is sufficiently smooth near the $t$-axis, then one can show 
\[
 d \mu(t) = |u(\vec{0},t)|^2 dt.
\]
The major part of this subsection is devoted to the first step of argument. 

\paragraph{Notations} We start by introducing a few notations
\begin{definition} Let $i,j\in \{1,2,\cdots,d\}$. We define
\begin{align*}
   \mathbf{L}_{ij} u & = \left(x_i\cdot\frac{\partial}{\partial x_j} - x_j \cdot \frac{\partial}{\partial x_i}\right) u;\\
   \mathbf{L}'_\pm u & = \left(\partial_t \pm \frac{x}{|x|}\cdot \nabla\right) u = (\partial_t \pm \partial_r)u.
\end{align*}
\end{definition}
\begin{remark} \label{convolution1}
A straightforward calculation shows that the operators $\mathbf{L}_{ij}$ satisfy
\begin{itemize}
 \item[(i)] $\mathbf{L}_{ij}$ are commutative with operators $\mathbf{L}'_\pm$ and the multiplication of a radial function $f(|x|)$. 
 \begin{align*}
  &\mathbf{L}_{ij} f(|x|) = f(|x|) \mathbf{L}_{ij};& &\mathbf{L}_{ij} \mathbf{L}'_\pm = \mathbf{L}'_\pm \mathbf{L}_{ij}.&
 \end{align*}
 \item[(ii)] We also have the following identities 
 \begin{align*}
  &\sum_{1\leq i<j\leq d} \mathbf{L}_{ij}^2 u = \Delta u - \sum_{i,j=1}^d x_i x_j u_{ij} - (d-1)\sum_{i=1}^d x_i u_i;\\
  & \sum_{1\leq i<j\leq d} \left|\mathbf{L}_{ij} u\right|^2  = |x|^2 |\slashed{\nabla} u|^2 \;\Rightarrow\; 
  \sum_{1\leq i<j\leq d} \left|\mathbf{L}_{ij} \left(|x|^{\beta} u\right)\right|^2 = |x|^{2+2\beta} |\slashed{\nabla} u|^2.
 \end{align*}
 \item[(iii)]  the following integral always vanishes for any differentiable function $g(x,t)$ and spatially radially symmetric region $\Omega$
 \[
  \iint_\Omega \mathbf{L}_{ij} g(x,t) dx dt = \int_{\partial \Omega} \left(x_i g(x,t) \vec{e}_j - x_j g(x,t) \vec{e}_i\right)\cdot d\mathbf{S} = 0.
 \]
 Because the normal vector of $\partial \Omega$ is alway orthogonal to $x_i \vec{e}_j - x_j \vec{e}_i$. 
\end{itemize}
\end{remark}

\paragraph{Divergence of vector field} Before we may apply Gauss' formula, we need to calculate the divergence
\begin{align*}
 \hbox{div}\mathbf{V}_- & = \left(\frac{d-1}{|x|}\!+\!\frac{x}{|x|}\cdot \nabla\right)\left[-\frac{\left|\mathbf{L}_+ u\right|^2}{4} + \frac{1}{2} e'(x,t)\right] + \frac{\partial}{\partial t}\left[\frac{\left|\mathbf{L}_+ u\right|^2}{4} + \frac{1}{2} e'(x,t)\right]\\
 & = \left(\partial_t \!+\! \frac{d-1}{|x|}\!+\!\frac{x}{|x|}\cdot \nabla\right) \frac{e'(x,t)}{2} + \left(\partial_t \!-\! \frac{d-1}{|x|}\!-\! \frac{x}{|x|}\cdot \nabla\right)\frac{|\mathbf{L}_+ u|^2}{4}\\
 & = J_1 + J_2.
\end{align*}
\noindent We first recall the definitions of the operators $\mathbf{L}'_\pm$, $\mathbf{L}_\pm$ and observe the following identities about operator compositions. These will be used frequently in our argument. 
\begin{align}
 &\partial_t \pm \left(\frac{\beta}{|x|} + \frac{x}{|x|}\cdot \nabla \right) = |x|^{-\beta} \mathbf{L}'_{\pm} |x|^\beta; & &\mathbf{L}_\pm = \pm |x|^{-\frac{d-1}{2}}\mathbf{L}'_\pm |x|^{\frac{d-1}{2}}.& \label{commutator1}
\end{align}
\paragraph{The term $J_1$} We may use the operator composition identities above to calculate $J_1$. Here we use the notation $w = |x|^{(d-1)/2} u$ and the identity $\displaystyle |x|^{d+1} |\slashed{\nabla} u|^2 =  \sum_{1\leq i<j\leq d} |\mathbf{L}_{ij} w|^2$ given in part (ii) of Remark \ref{convolution1}.
\begin{align*}
 J_1 & = \left(\partial_t \!+\! \frac{d-1}{|x|}\!+\!\frac{x}{|x|}\cdot \nabla\right)\left(\frac{1}{4}|\slashed{\nabla} u|^2 + \frac{\lambda_d}{4}\cdot \frac{|u|^2}{|x|^2} + \frac{1}{2(p+1)}|u|^{p+1} \right)\\
 & = \left(\partial_t \!+\! \frac{d+1}{|x|}\!+\!\frac{x}{|x|}\cdot \nabla\right)\frac{|\slashed{\nabla} u|^2}{4} - \frac{1}{2}\cdot \frac{|\slashed{\nabla} u|^2}{|x|} \\
 & \qquad \qquad + |x|^{-(d-1)}\mathbf{L}'_+ \left[|x|^{(d-1)}\left(\frac{\lambda_d}{4}\cdot \frac{|u|^2}{|x|^2} + \frac{1}{2(p+1)}|u|^{p+1}\right)\right]\\
 & = \frac{|x|^{-(d+1)}}{4} \mathbf{L}'_+ \left(|x|^{d+1} |\slashed{\nabla} u|^2\right) - \frac{1}{2}\cdot \frac{|\slashed{\nabla} u|^2}{|x|}\\
 & \qquad \qquad + |x|^{-(d-1)}\mathbf{L}'_+ \left(\frac{\lambda_d}{4}\cdot \frac{|w|^2}{|x|^2} + \frac{1}{2(p+1)}\cdot\frac{|w|^{p+1}}{|x|^{(d-1)(p-1)/2}} \right)\\
 & = \frac{|x|^{-(d+1)}}{4} \mathbf{L}'_+ \left(\sum_{1\leq i<j\leq d} |\mathbf{L}_{ij} w|^2\right) + J_3 - M(x,t)\\
 & = \frac{|x|^{-(d+1)}}{2} \sum_{1\leq i<j\leq d} (\mathbf{L}_{ij} w)(\mathbf{L}_{ij} \mathbf{L}'_+ w) + J_3 - M(x,t).
\end{align*}
The notation $J_3$ above represents
\[
 J_{3} = \frac{\lambda_d}{4}|x|^{-(d+1)} \mathbf{L}'_+(|w|^2)  + \frac{1}{2(p+1)}|x|^{-\frac{(d-1)(p+1)}{2}} \mathbf{L}'_+ (|w|^{p+1}).
\]
\paragraph{The term $J_2$} We can also use identities \eqref{commutator1} to calculate $J_2$
\begin{align*}
 J_2  = \frac{|x|^{-(d-1)}}{4} \mathbf{L}'_- \left| |x|^{\frac{d-1}{2}} \mathbf{L}_+ u\right|^2 
 = \frac{|x|^{-(d-1)}}{4} \mathbf{L}'_-\left|\mathbf{L}'_+ w\right|^2
 = \frac{|x|^{-(d-1)}}{2} \left[\mathbf{L}'_+ w\right] \cdot \left[\mathbf{L}'_- \mathbf{L}'_+ w\right]
\end{align*}
A straight forward calculation shows 
\begin{align*}
 \mathbf{L}'_- \mathbf{L}'_+ w  = |x|^{\frac{d-5}{2}}\left[|x|^2 u_{tt} - \lambda_d u - (d-1)\sum_{i=1}^d x_i u_i - \sum_{i,j=1}^d x_i x_j u_{ij}\right].
\end{align*}
We may plug in the equation $u_{tt} = \Delta u - |u|^{p-1}u$, use the first identity in part (ii) of Remark \ref{convolution1} and obtain
\begin{align*}
 \mathbf{L}'_- \mathbf{L}'_+ w & = |x|^{\frac{d-5}{2}}\left[\sum_{1\leq i<j\leq d} \mathbf{L}_{ij}^2 u -|x|^2 |u|^{p-1}u - \lambda_d u \right]\\
 & = |x|^{-2} \sum_{1\leq i<j\leq d} \mathbf{L}_{ij}^2 w - |x|^{-\frac{(p-1)(d-1)}{2}} |w|^{p-1}w - \lambda_d |x|^{-2} w
\end{align*}
We plug this in the expression of $J_2$ and obtain
\begin{align*}
 J_2 & = \frac{|x|^{-(d+1)}}{2}(\mathbf{L}'_+ w)\sum_{1\leq i<j\leq d} \mathbf{L}_{ij}^2 w - J_{3}.
\end{align*}
Now we combine $J_1$ with $J_2$ 
\begin{align}
 \hbox{div} \mathbf{V}_- & = \sum_{1\leq i<j\leq d} \mathbf{L}_{ij} \left[\frac{|x|^{-(d+1)}}{2} (\mathbf{L}_{ij} w) (\mathbf{L}'_+ w)\right] - M(x,t). \label{div formula}
\end{align}
\paragraph{Gauss' Formula} Now we apply Gauss' formula
\begin{align*}
 \int_{\partial \Omega} \mathbf{V}_{-} \cdot d\mathbf{S} & = \iint_{\Omega} \hbox{div} \mathbf{V}_- dxdt \\
 & = \iint_{\Omega} \left\{\sum_{1\leq i<j\leq d} \mathbf{L}_{ij} \left[\frac{|x|^{-(d+1)}}{2} (\mathbf{L}_{ij} w) (\mathbf{L}'_+ w)\right] - M(x,t)\right\} dx dt\\
 & = - \iint_\Omega M(x,t) dx dt.
\end{align*}
The integral of $\mathbf{L}_{ij}\left[\cdot\right]$ over a radially symmetric region vanishes, according to Remark \ref{convolution1}. This finishes the proof for inward energy. The outward energy can be dealt with in the same manner.
\begin{remark} \label{smooth approximation flux formula}
 Again we assume that the solution $u$ is sufficiently smooth in the proof above. When necessary we may apply the same smooth approximation techniques as in the proof of Morawetz estimates. The error term is 
 \begin{align*}
  Err = \frac{1}{2} \iint_\Omega |x|^{-(d-1)} \left(\mathbf{L}'_+ w_\eps\right)  |x|^{\frac{d-1}{2}} (F_\eps + |u_\eps|^{p-1} u_\eps) dx dt 
  = \frac{1}{2} \iint_\Omega (\mathbf{L}_+ u_\eps) (F_\eps + |u_\eps|^{p-1} u_\eps) dx dt,
 \end{align*}
 which is dominated by
 \[
  Err \lesssim_1 \|\mathbf{L}_+ u_\eps\|_{L^\infty L^2} \|F_\eps + |u_\eps|^{p-1} u_\eps\|_{L^1 L^2} \lesssim_1 E^{1/2}  \|F_\eps + |u_\eps|^{p-1} u_\eps\|_{L^1 L^2}.
 \]
 This clearly vanishes as $\eps \rightarrow 0^+$.
\end{remark}
\begin{remark} \label{nonradial regions}
 One may also consider energy flux formula of inward/outward energy in a region that is not necessarily radially symmetric. In this case we need to substitute the original vector fields $\mathbf{V}_-$ by $\mathbf{V}_- - \frac{1}{2}(\mathbf{L}_+ u)\slashed{\nabla} u$. In fact we have
 \[
  \hbox{div} \left(\mathbf{V}_- - \frac{1}{2} (\mathbf{L}_+ u)\slashed{\nabla} u \right) = -M(x,t).
 \]
 In order to show this we combine identity \eqref{div formula} with 
 \begin{align*}
  \hbox{div} \left[(\mathbf{L}_+ u)\slashed{\nabla} u\right] & = \hbox{div} \left[(\mathbf{L}_+ u)\nabla u - (\mathbf{L}_+ u) \left(\frac{x}{|x|}\cdot \nabla u\right) \frac{x}{|x|}\right] \\
  & = \hbox{div} \left\{\sum_{j=1}^d \left[\sum_{i=1}^d \frac{x_i}{|x|} \cdot \frac{x_i u_j - x_j u_i}{|x|} \cdot (\mathbf{L}_+ u)\right] \vec{e}_j\right\}\\
  & = \hbox{div} \left\{\sum_{j=1}^d \left[\sum_{i=1}^d x_i |x|^{-(d+1)} (\mathbf{L}_{ij} w) (\mathbf{L}'_+ w) \right] \vec{e}_j\right\}\\
  & = \sum_{1\leq i<j\leq d} \mathbf{L}_{ij} \left[ |x|^{-(d+1)}(\mathbf{L}_{ij} w) (\mathbf{L}'_+ w)\right].
 \end{align*}
 Similarly we substitute $\mathbf{V}_+$ by $\mathbf{V}_+ + \frac{1}{2}(\mathbf{L}_- u) \slashed{\nabla} u$ in the case of outward energy. 
 \end{remark}
 
\subsection{Asymptotic behaviour of inward/outward energy}

As in the 3-dimensional case, we may prove the following propositions regarding the space-time distribution of energy thus prove Theorem \ref{main 1}. The key observation is that the Morawetz density function $M(x,t)$ and the non-directional component of energy density $e'(x,t)$ are comparable to each other
\[
 e'(x,t) \simeq |x|\cdot M(x,t).
\]
We omit the details of proof since the argument is exactly the same as in the 3-dimensional case. Please see \cite{shen3dnonradial} for the details. In all of the following propositions the pair $(p,d)$ is always assumed to satisfy (A1).
\begin{proposition}[Monotonicity and limit of inward/outward energies] 
Let $u$ be a solution to (CP1) with a finite energy. The inward energy $E_-(t)$ is a decreasing function of $t$; the outward energy $E_+(t)$ is an increasing function of $t$. In addition we have the following limits
 \begin{align*}
  &\lim_{t \rightarrow \pm \infty} E_\pm (t) = E;& &\lim_{t \rightarrow \mp \infty} E_{\pm} (t) = 0.&\\
  &\lim_{t \rightarrow \pm \infty} \int_{\Rm^d} \left(|\slashed{\nabla} u|^2 + |u|^{p+1} + \lambda_d \cdot \frac{|u|^2}{|x|^2}\right) dx = 0.& &&
 \end{align*}
\end{proposition}
\begin{proposition} \label{rediscover of Morawetz}
We have an expression of inward energy at time $t_0$ in terms of the measure $\mu$, as introduced in the energy flux formula, and the Morawetz integral
\[
 E_-(t_0) = c_d \mu([t_0,\infty)) + \int_{t_0}^\infty \int_{\Rm^d} M(x,t) dx dt.
\]
We may also rediscover the Morawetz estimates
\[
 c_d \mu(\Rm) + \iint_{\Rm^d \times \Rm} M(x,t) dx dt = E. 
\]
\end{proposition}
\begin{proposition}[Asymptotic behaviour of fluxes through light cones] Let $Q_-^- (s)$ and $Q_+^+(\tau)$ be energy fluxes through light cones, as in Definition \ref{def fluxes light cone}. Then the following limits hold  
\begin{align*}
 &\lim_{s \rightarrow +\infty} Q_-^-(s) = 0;& &\lim_{\tau \rightarrow -\infty} Q_+^+(\tau) = 0.&
\end{align*}
\end{proposition}
\begin{proposition}[Travelling speed of energy]
 We have the limits
 \[
  \lim_{t \rightarrow \pm \infty} E_\pm (t;\{x\in \Rm^d: |x|<c|t|\}) = 0
\]
for any constant $c \in (0,1)$. It immediately follows that
 \begin{align*}
  \lim_{t\rightarrow \pm \infty} \int_{|x|<c|t|} \left(\frac{1}{2}|\nabla u(x,t)|^2 + \frac{1}{2}|u_t(x,t)|^2 + \frac{1}{p+1}|u(x,t)|^{p+1} \right) dx = 0.
 \end{align*}
\end{proposition}

\subsection{Weighted Morawetz estimates}
As an application of our inward/outward energy theory, we may prove weighted Morawetz estimates and the decay estimates of inward/outward energy, as given in part (a) of Theorem \ref{main 2}. We start by 
\begin{proposition}[General Weighted Morawetz] \label{Weighted Morawetz General}
Assume that $(d,p)$ satisfies (A1) and $0<\gamma<1$.  Let $a(r) \in C([0,\infty))$ be a function satisfying
\begin{itemize}
 \item $a(r)$ is absolutely continuous in $[0,R]$ for any $R \in \Rm^+$;
 \item Its derivative satisfies $0 \leq a'(r) \leq \gamma a(r)/r$ almost everywhere $r>0$.
\end{itemize}
 If $u$ is a solution to (CP1) with a finite energy so that 
 \begin{align*}
  K_1 =  \int_{\Rm^d} a(|x|) \left[\frac{\left|(\mathbf{L}_+ (u_0,u_1))(x)\right|^2}{4} + \frac{\lambda_d |u_0(x)|^2}{4|x|^2} + \frac{|\slashed{\nabla} u_0(x)|^2}{4} +\frac{|u_0(x)|^{p+1}}{2(p+1)} \right] dx < \infty,
 \end{align*}
 then we have 
 \begin{align*}
   \int_0^\infty a(t) d\mu(t) + \iint_{\Rm^d \times \Rm^+} \!\! a(|x|+t)\left(\frac{|\slashed{\nabla} u(x,t)|^2}{|x|} + \frac{\lambda_d |u(x,t)|^2}{|x|^3} + \frac{|u(x,t)|^{p+1}}{|x|} \right) dxdt \lesssim_{d,p,\gamma} K_1,
 \end{align*}
 where $\mu$ is the measure introduced in the energy flux formula.
\end{proposition} 

\begin{remark} \label{gamma 1 remark}
Let $\gamma = 1$. If $(p,d)$, $a(r)$ and $u$ still satisfy conditions in Proposition \ref{Weighted Morawetz General}, then we also have a weighted Morawetz estimate
 \begin{align*}
   \int_0^\infty a(t) d\mu(t) + \iint_{\Rm^d \times \Rm^+} \!\! \frac{t\cdot a(|x|+t)}{|x|+t}\left(\frac{|\slashed{\nabla} u(x,t)|^2}{|x|} + \frac{\lambda_d |u(x,t)|^2}{|x|^3} + \frac{|u(x,t)|^{p+1}}{|x|} \right) dxdt \lesssim_{d,p} K_1.
 \end{align*}
\end{remark}

\noindent Proposition \ref{Weighted Morawetz General} and Remark \ref{gamma 1 remark} can be proved in the same way as in 3-dimensional case. We will not repeat this proof in this work. Readers may see \cite{shen3dnonradial} for the details of the proof. 

\paragraph{Proof of Part (a), Theorem \ref{main 2}} Now we are ready to prove the weighted Morawetz estimates given in Theorem \ref{main 2}. We apply Proposition \ref{Weighted Morawetz General} with $a(r) = r^\kappa$ and $\gamma = \kappa$. According to Remark \ref{relationship of u w energy}, we have 
\begin{align*}
 K_1 &\leq  \int_{\Rm^d} |x|^\kappa \left[\frac{\left|\mathbf{L}_+ (u_0,u_1)\right|^2}{4} + \frac{\left|\mathbf{L}_- (u_0,u_1)\right|^2}{4}  + \frac{\lambda_d |u_0|^2}{2|x|^2} + \frac{|\slashed{\nabla} u_0|^2}{2} + \frac{|u_0|^{p+1}}{p+1} \right] dx \\
 & \leq  \int_{\Rm^d} |x|^{\kappa}\left[\frac{1}{2}|\nabla u_0|^2 + \frac{1}{2}|u_1|^2 + \frac{1}{p+1}|u_0|^{p+1}\right]dx \leq E_\kappa(u_0,u_1)< + \infty.
\end{align*}
The conclusion of Proposition \ref{Weighted Morawetz General} is  
\[
 \int_0^\infty t^\kappa d\mu(t) + \int_{0}^\infty \int_{\Rm^d} (t+|x|)^\kappa \left(\frac{|\slashed{\nabla} u|^2}{|x|} + \lambda_d \frac{|u|^2}{|x|^3} + \frac{|u|^{p+1}}{|x|}\right) dx dt \lesssim_{d,p,\kappa} K_1 \leq E_\kappa(u_0,u_1).
\]
This immediately gives the weighted Morawetz estimate for positive time $t>0$. The negative time direction can be dealt with in the same manner since the wave equation is time-reversible. Next we prove the decay estimate of inward/outward energy. Without loss of generality let us prove the estimates of $E_-(t)$ for positive time $t$. We first write $E_-(t)$ in terms of $\mu$ and Morawetz integral by Proposition \ref{rediscover of Morawetz}. 
\begin{equation}
 E_-(t)  = c_d \mu([t,\infty)) + \int_{t}^\infty \int_{\Rm^d} M(x,t') dx dt'.  \label{expression of E minus}
\end{equation}
This immediately gives an upper bound of $E_-(t)$
\begin{align*}
 E_-(t) & \leq t^{-\kappa} c_d \int_t^{\infty} |t'|^{\kappa} d\mu(t') + t^{-\kappa} \int_t^\infty \int_{\Rm^d} (t'+|x|)^{\kappa} M(x,t') dx dt'\\
 & \lesssim_{d,p,\kappa} t^{-\kappa} E_{\kappa} (u_0,u_1).
\end{align*}
Here we apply the weighted Morawetz estimate. Finally let us show $E_-(t) \in L^{1/\kappa} ([0,\infty))$. We simply find an upper bound of $E_-(t)$ by \eqref{expression of E minus}  
\[
 E_-(t) \leq c_d \int_t^\infty (t')^{-\kappa} \cdot (t')^\kappa d\mu(t') + \int_{t}^\infty (t')^{-\kappa}\left(\int_{\Rm^d} (t'+|x|)^\kappa M(x,t') dx\right) dt',
\]
then apply the following lemma with $d \mu'(t') = (t')^\kappa d\mu(t')$ and $d\mu'(t') = \left(\int_{\Rm^d} (t'+|x|)^\kappa M(x,t') dx \right) dt'$, respectively. 

\begin{lemma}
 Let $\mu'$ be a continuous, nonnegative, finite measure on $[0,\infty)$ and $\kappa \in (0,1)$ be a constant. Then the function $f: \Rm^+ \rightarrow \Rm$ defined by $f(x) = \int_x^\infty y^{-\kappa} d\mu'(y)$ satisfies $f \in L^{1/\kappa} (\Rm^+)$. In fact we have $\|f\|_{L^{1/\kappa}(\Rm^+)} \leq \mu'(\Rm^+)$.
\end{lemma}
\begin{proof}
 We start by finding an upper bound of $f(x)$
 \begin{align*}
  f(x)  =  \int_x^\infty (1\cdot y^{-\kappa}) d\mu'(y) & \leq  \left(\int_x^\infty 1^{1/(1-\kappa)} d\mu'(y)\right)^{1-\kappa} \left(\int_x^\infty \left(y^{-\kappa}\right)^{1/\kappa} d\mu'(y)\right)^{\kappa} \\
  & \leq \mu'(\Rm^+)^{1-\kappa} \left(\int_x^\infty y^{-1} d\mu'(y)\right)^{\kappa}
 \end{align*}
 Therefore we have
 \begin{align*}
  \int_0^\infty |f(x)|^{1/\kappa} dx & \leq \mu'(\Rm^+)^{(1-\kappa)/\kappa} \int_0^\infty \left(\int_x^\infty y^{-1} d\mu'(y)\right) dx \\
  & =  \mu'(\Rm^+)^{(1-\kappa)/\kappa} \int_0^\infty \left(\int_0^y y^{-1} dx \right) d\mu'(y)\\
  & = \mu'(\Rm^+)^{1/\kappa}.
 \end{align*}
\end{proof}
\section{Scattering Theory}

\begin{lemma} \label{abstract inter}
Assume that $(d,p)$ satisfies (A1). Let $q_1,r_1,q_2,r_2,k_1,k_2$ be constants so that
\begin{itemize}
 \item The pair $(q_2,r_2)$ is $1$-admissible;
 \item $1<q_1,r_1<\infty$, $k_1,k_2>0$ and the following identities hold
 \begin{align*}
  &k_1+k_2 = p-1;& &k_1/q_1 +k_2/q_2 = 2/(d+1);& &k_1/r_1+k_2/r_2 = 2/(d+1).& 
 \end{align*}
\end{itemize} 
 If $u$ is a finite-energy solution to (CP1) with $\|u\|_{L^{q_1} L^{r_1}(\Rm^+ \times \Rm^d)} < + \infty$, then we also have
 \[
    \|u\|_{S_{d,p}(\Rm)} \doteq \|u\|_{L^{(d+1)(p-1)/2} L^{(d+1)(p-1)/2} (\Rm^+ \times \Rm^d)}< +\infty.
 \]
In addition, the solution scatters in the energy space as $t \rightarrow +\infty$. More precisely, there exists $(v_0^+,v_1^+) \in \dot{H}^1 \times L^2(\Rm^d)$, so that 
\[
 \lim_{t \rightarrow +\infty} \left\|\begin{pmatrix}u(\cdot,t)\\ \partial_t u(\cdot,t)\end{pmatrix} - \mathbf{S}_L (t) \begin{pmatrix} v_0^+\\ v_1^+\end{pmatrix}\right\|_{\dot{H}^{1} \times L^2} = 0.
\]
\end{lemma}  
\begin{proof}
 We may follow the same argument as in the 3-dimensional case. Thus we only give a sketch of proof. Please see \cite{shen3dnonradial} for more details. We first define a function $g_{t_0}$ on $[t_0,\infty)$ by
 \[
 g_{t_0}(T) = \|u\|_{L^{q_2} L^{r_2}([t_0,T]\times \Rm^d)} +  \|D_x^{1/2} u\|_{L^{\frac{2(d+1)}{d-1}} L^{\frac{2(d+1)}{d-1}} ([t_0,T]\times \Rm^3)}.
\]
The value of $g_{t_0}(T)$ is always finite for all $0\leq t_0<T<\infty$ by the Strichartz estimates. We may combine the Strichartz estimates, the fractional chain rule and the energy conservation law to obtain
 \begin{align*}
  g_{t_0}(T) & \lesssim \|(u(\cdot,t_0),u_t(\cdot,t_0))\|_{\dot{H}^1 \times L^2} + \|D_x^{1/2} \left(-|u|^{p-1} u\right)\|_{L^\frac{2(d+1)}{d+3} L^\frac{2(d+1)}{d+3} ([t_0,T]\times \Rm^d)}\\
  & \lesssim E^{1/2} + \|D_x^{1/2} u\|_{L^{\frac{2(d+1)}{d-1}} L^{\frac{2(d+1)}{d-1}}([t_0,T]\times \Rm^d)} \|u\|_{L^{p_1}L^{r_1}([t_0,T]\times \Rm^d)}^{k_1} \|u\|_{L^{q_2} L^{r_2}([t_0,T]\times \Rm^d)}^{k_2} \\
  & \lesssim E^{1/2} + \|u\|_{L^{q_1}L^{r_1}([t_0,T]\times \Rm^d)}^{k_1} \left(g_{t_0}(T)\right)^{k_2+1}.
 \end{align*}
 In other words, there exists a constant $C$ independent of $t_0, T$, so that
 \begin{equation} \label{self inequality}
  g_{t_0}(T) \leq CE^{1/2} + C\|u\|_{L^{p_1}L^{r_1}([t_0,T]\times \Rm^d)}^{k_1} \left(g_{t_0}(T)\right)^{k_2+1}.
 \end{equation}
 Please note that $\|u\|_{L^{q_1} L^{r_1}(\Rm^+ \times \Rm^d)} < +\infty$ implies $\displaystyle \lim_{t_0\rightarrow +\infty}  \|u\|_{L^{q_1} L^{r_1}([t_0,\infty) \times \Rm^d)} = 0$. When $t_0$ is sufficiently large, we may apply a continuity argument to conclude $g_{t_0}(T) < 2CE^{1/2}$ for all $T>t_0$. As a result we have 
 \begin{align*}
  \|u\|_{L^{q_2} L^{r_2}(\Rm^+\times \Rm^d)}+\|D_x^{1/2} u\|_{L^{\frac{2(d+1)}{d-1}} L^{\frac{2(d+1)}{d-1}} (\Rm^+\times \Rm^d)} & < +\infty; \\
  \Rightarrow  \|u\|_{L^{(d+1)(p-1)/2} L^{(d+1)(p-1)/2} (\Rm^+ \times \Rm^d)} & < +\infty.
 \end{align*}
 Finally we may verify the scattering of solution as $t \rightarrow +\infty$ by the Strichartz estimates, the completeness of the space $\dot{H}^1\times L^2$ and the fact that $\mathbf{S}_L(t)$ is a unitary operator in $\dot{H}^1 \times L^2$. 
\end{proof}

\paragraph{Proof of Theorem \ref{main 2} part (b)} The idea is to apply Lemma \ref{abstract inter} with the following constants $q_1,r_1,q_2,r_2, k_1,k_2$.
\begin{align*}
 &(q_1,r_1) = \left(\frac{p+1}{\kappa_0(d,p)}, p+1\right);& &(q_2,r_2) = \left(2, \frac{2d}{d-3}\right);& \\
 &k_1 = \frac{\frac{4d}{d+1} - (d-3)(p-1)}{\frac{2d}{p+1} -(d-3)};& &k_2 = \frac{-\frac{4d}{d+1} + 2d \cdot \frac{p-1}{p+1}}{\frac{2d}{p+1}-(d-3)}.&
\end{align*}
Here $\kappa_0(d,p) = \frac{(d+2)(d+3) - (d+3)(d-2)p}{(d-1)(d+3) - (d+1)(d-3)p}$ is defined in Theorem \ref{main 2}. We have already known
\[
 u \in L^\infty L^{p+1}(\Rm \times \Rm^d) \cap L^{(p+1)/\kappa} L^{p+1}(\Rm \times \Rm^d),
\]
by the energy conservation law and the decay estimates given in Part (a) of Theorem \ref{main 2}. Since we have $(p+1)/\kappa \leq (p+1)/\kappa_0(d,p) = q_1 < +\infty$ and $r_1=p+1$, we obtain $u \in L^{q_1} L^{r_1}(\Rm \times \Rm^d)$ by an interpolation. Now we can apply Lemma \ref{abstract inter} to conclude that $u$ scatters in the energy space as $t \rightarrow +\infty$ and $u\in S_{d,p}(\Rm^+)$. According to the scattering criterion Proposition \ref{scattering criterion}, the latter also implies that $u$ scatters in the critical Sobolev space $\dot{H}^{s_p} \times \dot{H}^{s_p-1}$ in the positive time direction. Finally an interpolation between these two spaces shows that the scattering happens in all spaces $\dot{H}^s \times \dot{H}^{s-1}$ for $s\in [s_p,1]$. The negative time direction can be dealt with in the same way since the wave equation is time-reversible. 

\begin{figure}[h]
 \centering
 \includegraphics[scale=1.0]{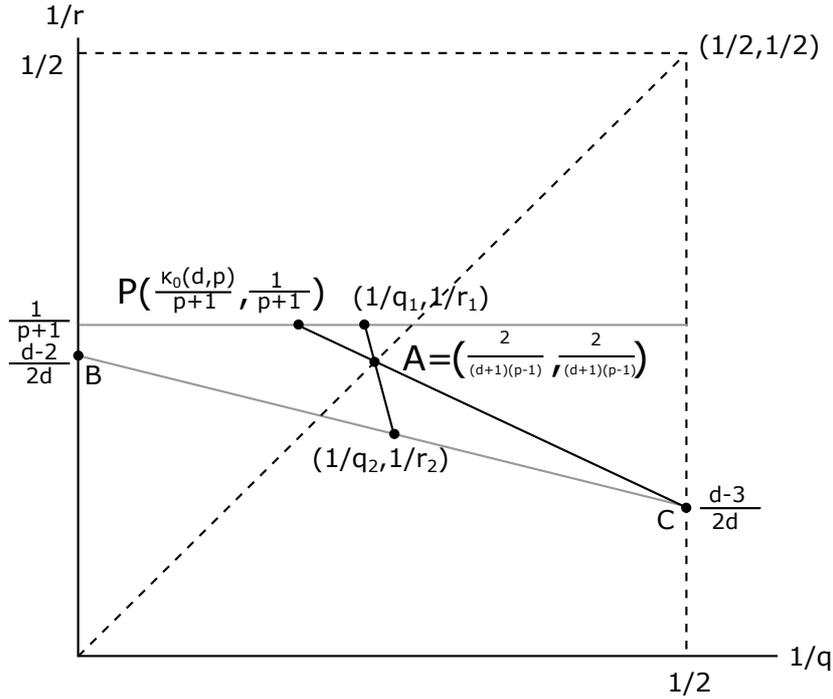}
 \caption{Best choices of admissible pair} \label{figure bestqr}
\end{figure}

\begin{remark}
The pair $(q_2,r_2) = (2,\frac{2d}{d-3})$ is our best choice when we prove the scattering theory by a combination of Lemma \ref{abstract inter} and a decay estimate $\|u\|_{L^{q_1} L^{p+1}(\Rm \times \Rm^d)} < +\infty$. This fact is visually displayed in figure \ref{figure bestqr}. In fact we are solving an optimization problem to minimize $1/q_1 = \kappa/(p+1)$ with constraints
\begin{itemize}
 \item The point $(1/q_1,1/r_1)$ is on the horizontal line $1/r_1 = 1/(p+1)$;
 \item The point $(1/q_2,1/r_2)$ is on the line segment $BC$, because of the $1$-admissible assumption;
 \item The point $A(\frac{2}{(d+1)(p-1)}, \frac{2}{(d+1)(p-1)})$ is on the line segment connecting $(1/q_1,1/r_1)$ and $(1/q_2,1/r_2)$.
\end{itemize} 
Simple observation makes clear the optimal choice:
\begin{align*}
 &(1/p_1,1/r_1) = \left(\frac{\kappa_0(d,p)}{p+1}, \frac{1}{p+1}\right),& &(1/q_2,1/r_2) = \left(\frac{1}{2}, \frac{d-3}{2d}\right).&
\end{align*}
\end{remark} 

\end{document}